\documentclass[10pt]{amsart}
\usepackage{amssymb}
\usepackage{amsbsy}
\usepackage{amscd}
\usepackage{amsmath}
\usepackage{amsthm}

\theoremstyle{plain}
 \newtheorem{thm}{Theorem}[section]
 \newtheorem{lem}[thm]{Lemma}
 \newtheorem{prop}[thm]{Proposition}
 \newtheorem{cor}[thm]{Corollary}
 
\theoremstyle{definition}
 \newtheorem{defn}{Definition}[section]
\theoremstyle{remark}
 \newtheorem{rem}{Remark}[section]
 \newtheorem{ex}{Example}[section]
 \newtheorem{claim}{Claim}[section]
\def\Bbb{\mathbb}
\def\frak{\mathfrak}
\def\cal{\mathcal}

\newcommand{ \Supp}{\operatorname{Supp}}
\newcommand{\Ext}{\operatorname{Ext}}
\newcommand{\Hom}{\operatorname{Hom}}

\newcommand{\im}{\operatorname{im}}

\newcommand{\rk}{\operatorname{rk}}

\newcommand{\NS}{\operatorname{NS}}
\newcommand{\coker}{\operatorname{coker}}
\newcommand{\Pic}{\operatorname{Pic}}
\newcommand{\ch}{\operatorname{ch}}
\newcommand{\td}{\operatorname{td}}

\newcommand{\Hilb}{\operatorname{Hilb}}
\newcommand{\Quot}{\operatorname{Quot}}

\newcommand{\Spec}{\operatorname{Spec}}

\newcommand{\WIT}{\operatorname{WIT}}

\newcommand{\Div}{\operatorname{Div}}

\font\b=cmr10 scaled \magstep5
\def\bigzerou{\smash{\lower1.7ex\hbox{\b 0}}}

\numberwithin{equation}{section}
\begin{document}

\title{
Twisted stability and Fourier-Mukai transform
}
\author{K\={o}ta Yoshioka}
 
\address{
Department of mathematics, Faculty of Science, Kobe University,
Kobe, 657, Japan}

\email{yoshioka@math.kobe-u.ac.jp}
 \maketitle

\section{Introduction}
Let $X$ be an abelian surface or a K3 surface over ${\Bbb C}$.
Mukai introduced a lattice structure $\langle \quad,\quad \rangle$
on 
$H^{ev}(X,{\Bbb Z}):=\oplus H^{2i}(X,{\Bbb Z})$
by 
\begin{equation}
\begin{split}
\langle x,y \rangle:=&-\int_X x^{\vee} \wedge y\\
=& \int_X(x_1 \wedge y_1-x_0 \wedge y_2-x_2 \wedge y_0),
\end{split} 
\end{equation}
where $x_i \in H^{2i}(X,{\Bbb Z})$ (resp. $y_i \in H^{2i}(X,{\Bbb Z})$)
is the $2i$-th component of $x$ (resp. $y$)
and $x^{\vee}=x_0-x_1+x_2$.
It is now called Mukai lattice.
For a coherent sheaf $E$ on $X$,
we can attach an element of $H^{ev}(X,{\Bbb Z})$
called Mukai vector $v(E):=\ch(E)\sqrt{\td_X}$,
where $\ch(E)$ is the Chern character of $E$ and 
$\td_X$ is the Todd class of $X$.
For a Mukai vector $v \in H^{ev}(X,{\Bbb Z})$ and an ample
divisor $H$, 
let $M_H(v)$ be the moduli space of stable sheaves $E$ of 
Mukai vector $v(E)=v$ and
$\overline{M}_H(v)$ the moduli space of semi-stable sheaves.
An ample divisor $H$ is general with respect to $v$, if the following
condition holds:
\begin{itemize}
\item[$(\natural)$]
for every $\mu$-semi-stable sheaf $E$ of $v(E)=v$,
if $F \subset E$ satisfies
$(c_1(F),H)/\rk F=(c_1(E),H)/\rk E$, then
$c_1(F)/\rk F=c_1(E)/\rk E$.
\end{itemize}

The preservation of stability
by Fourier-Mukai transform on $X$ was investigated by
many people (e.g. \cite{BBH:2}, \cite{B-M:1}, \cite{Mu:5},
\cite{Y:7}).
In \cite{Y:7},
we introduced twisted degree of coherent sheaf $E$
by $\deg_G(E)=\deg(E \otimes G^{\vee})$,
where $G$ is a vector bundle on $X$.
Then we showed that Fourier-Mukai transform preserves Gieseker semi-stability, 
if twisted degree is 0 and the polarization $H$ is general.
In this paper, we shall generalize our results to
 the case where $H$ is not general.
In this case, Fourier-Mukai transform does not preserve 
Gieseker semi-stability.
This fact is closely related to the following fact:
If $H$ is not general, than Gieseker semi-stability is not
preserved by the twisting $E \mapsto E \otimes L$, where
$L$ is a line bundle. 
Thus Gieseker semi-stability depends on the choice of $L$.
In order to understand this phenomenon,
Matsuki and Wentworth \cite{M-W:1}(also by Ellingsrud and
G\"{o}ttsche \cite{E-G:1} and
Friedman and Qin \cite{F-Q:1})
introduced $L$-twisted semi-stability, where $L$ is a ${\Bbb Q}$-line bundle.
Hence we shall propose a formulation for our problem
by using twisted semi-stability.
In section \ref{sect:K3}, we shall show that Fourier-Mukai transform preserves
suitable twisted semi-stability, if $X$ is an abelian surface
(Theorem \ref{thm:00/08/17}).

In \cite{Y:7}, we showed that $M_H(v)$ is
deformation equivalent to a moduli space of torsion free sheaves of 
rank 1, if $v$ is primitive and the polarization $H$ is general.
In section \ref{sect:deform}, 
we shall give another proof of this result
by using results proved in section \ref{sect:K3}.
Moreover we shall show the following.
\begin{thm}
Let $X$ be an abelian surface or a K3 surface.
Let $v \in H^{ev}(X,{\Bbb Z})$ be a Mukai vector of $\rk v>0$.
Then 
$\overline{M}_H(v)$ is a normal variety, if $\langle v^2 \rangle>0$
and $H$ is general with respect to $v$.
\end{thm}

In section \ref{sect:enriques}, we shall consider Fourier-Mukai transform
on an Enriques surface
associated to $(-1)$-reflection.
In particular, we shall show a similar result to 
Theorem \ref{thm:00/08/17}
(Proposition \ref{prop:00/08/17}).
As an application, we shall compute Hodge polynomials of some
moduli spaces (Theorem \ref{thm:enriques}).
We also discuss a relation to Montonen-Olive duality in Physics
(cf. \cite{V-W:1}).

We are also interested in Fourier-Mukai transform on elliptic
surfaces.
Let $\pi:X \to C$ be an elliptic surface with a $0$-section.
Then a compactification of the relative Jacobian is isomorphic to $X$ and
there is a universal family ${\cal P}$ on $X \times_C X$.
We regard ${\cal P}$ as a sheaf on $X \times X$ and
consider Fourier-Mukai transform defined by ${\cal P}$.
Assume that every fiber is irreducible.
Then the preservation of stability was investigated in
\cite{Br:1}, \cite{H-M:1}, \cite{J-M:1}, \cite{Y:7}.
In particular, semi-stable sheaf of relative degree 0 maps to
a semi-stable sheaf of pure dimension 1.
In \cite{Br:1}, Bridgeland also treated Fourier-Mukai transform induced by
a relative moduli space of stable sheaves on fibers.
In order to generalize our result \cite[Thm. 3.15]{Y:7} to
this situation, we need to consider twisted semi-stability
for purely 1-dimensional sheaves.
In section \ref{sect:1-dim}, 
we introduce twisted stability for purely 1-dimensional
sheaves and give some properties which are similar to
results in \cite{Y:2}. 
Then we can show in Theorem \ref{thm:FM} that
a twisted semi-stable sheaf of twisted relative degree 0 maps to
a twisted semi-stable sheaf of pure dimension 1.
As an application, we can compute Hodge numbers of some moduli spaces.

\begin{thm}
Let $\pi:X \to C$ be an elliptic surface with a section $\sigma$.
Assume that every fiber is irreducible.
$f$ denotes a fiber of $\pi$.
Let $M_{\sigma+kf}(r,c_1,\chi)$ be the moduli space of stable sheaves $E$ of
$(\rk(E),c_1(E),\chi(E))=(r,c_1,\chi)$ with respect to $\sigma+kf$.
Then 
\begin{equation}
h^{p,q}(M_{\sigma+kf}(r,c_1,\chi))=h^{p,q}(\Pic^0(X) \times \Hilb_X^n),
\end{equation}
if $(r,(c_1,f))=1$ and $k \gg 0$,
where $2n+h^1({\cal O}_X)=\dim M_{\sigma+kf}(r,c_1,\chi)$.
\end{thm}
Indeed we introduced twisted semi-stability for purely 1-dimensional sheaves
to prove this theorem.

We also construct moduli spaces of twisted semi-stable sheaves
by using Ellingsrud and G\"{o}ttsche's method:
They used moduli spaces of parabolic semi-stable torsion free
sheaves constructed by Yokogawa \cite{Yk:1}.
Since his construction only works for
parabolic semi-stable torsion free sheaves on smooth projective scheme,
we need to construct moduli spaces of parabolic semi-stable sheaves
of pure dimension 1, 
which will be done in appendix (Theorem \ref{thm:parabolic}).

\section{Preliminaries}\label{sect:pre}

{\it Notation.}

Throughtout this note, we use the following notations.
Let $X$ be a smooth surface.
For a scheme $S$, we denote the projection
$S \times X \to S$ by
$p_S$.

Let ${\Bbb Q}[[1/x,1/y]]$ be the formal power series ring of
two variable $1/x,1/y$ and
$R$ the localization of ${\Bbb Q}[[1/x,1/y]]$ 
by $1/(xy)$.

\subsection{Virtual Hodge polynomial}

For a variety $Y$ over ${\Bbb C}$, cohomology with compact support
$H^*_c(Y,{\Bbb Q})$ has a natural mixed Hodge structure.
Let $e^{p,q}(Y):=\sum_k(-1)^k h^{p,q}(H_c^k(Y))$ be the virtual Hodge number
and $e(Y):=\sum_{p,q}e^{p,q}(Y)x^p y^q$
the virtual Hodge polynimial of $Y$.
For more details, see [Ch, 0.1].
By the properties of virtual Hodge polynomials,
we can extend the definition of $e(Y)$ to constructible sets.
 
Let $S$ be a bounded set of coherent sheaves on $X$ with
Hilbert polynomial $h(x)$.
Under some conditions, we shall define
the virtual Hodge polynomial of $S$ as an element of $R$.
So it is not a polynomial in general.
 Let $H$ be an ample divisor on $X$.
Since $S$ is bounded,
there is an integer $m$ such that
for every element $E$,
\begin{enumerate}
\item
$E(mH)$ is generated by global sections,
\item
$H^i(X,E(mH))=0$ for $i>0$.
\end{enumerate}
We set $N:=h(m)$.
We shall consider the quot scheme
$Q:=\Quot_{{\cal O}_X(-mH)^{\oplus N}/X/{\Bbb C}}^{h(x)}$.
Let ${\cal O}_{Q \times X}(-mH)^{\oplus N}
 \to {\cal Q}$ be the universal quotient.
We assume that 
$$
Q_S:=\{q \in Q|{\cal Q}_q \in S, \text{$H^0(X,{\cal O}_X^{\oplus N}) \to 
H^0(X,{\cal Q}_q(mH))$ is isomorphic} \}
$$ 
is a constructible set.
We shall show that this condition does not depend on the choice of $Q$.
Let $T$ be a scheme and ${\cal E}$ a coherent sheaf on $T \times X$
such that ${\cal E}$ is flat over $T$.
Let $T^0$ be the open subscheme of $T$ 
consisting of point $t \in T_S$ such that
\begin{enumerate}
\item
${\cal E}_t(mH)$ is generated by global sections,
\item
$H^i(X,{\cal E}_t(mH))=0$ for $i>0$.
\end{enumerate}
By base change theorem, 
$p_{T^0*}({\cal E}(mH))$ is a locally free sheaf on $T^0$.
Let $\tau:{\cal P} \to T^0$ be the associated principal 
$GL(N)$ bundle over $T^0$.
Then there is a surjective homomorphism
${\cal O}_{{\cal P} \times X}(-mH)^{\oplus N}
 \to (\tau \times id_X)^*{\cal E}$.
Hence we get a morphism $\eta:{\cal P} \to Q$.
We set $T_S:=\{t \in T^0|{\cal E}_t \in S \}$.
Then we see that $T_S=\tau(\eta^{-1}(Q_S))$.
Thus $T_S$ is constructible, if $Q_S$ is constructible.
In particular, constructibility does not depend on the choice of $Q$.
\begin{defn}\cite[Defn. 1.1]{Y:remark}
\begin{enumerate}
\item
$S$ is constructible, if $Q_S$ is constructible.
\item
For a constructible set $S$, we define the virtual Hodge polynomial by
$$
e(S):=\frac{e(Q_S)}{e(GL(N))} \in R.
$$ 
\end{enumerate}
\end{defn}
It is easy to see that this definition does not depend on the choice of $Q$.
Let ${\cal M}_S$ be the associated
substack of the stack of coherent sheaves.
We define the virtual Hodge polynomial $e({\cal M}_S)$ by
$e(S)$. 

\subsection{Fourier-Mukai transform}
Let $K(X)$ be the Grothendieck group of $X$.
For $x \in K(X)$, we set 
\begin{equation}
\gamma(x):=(\rk x,c_1(x),\chi(x)) \in {\Bbb Z} \oplus \NS(X) \oplus
{\Bbb Z}.
\end{equation}
Then $\gamma:K(X) \to {\Bbb Z} \oplus \NS(X) \oplus
{\Bbb Z}$ is a surjective homomorphism and $\ker \gamma$ is generated
by ${\cal O}_X(D)-{\cal O}_X$ and
${\Bbb C}_P-{\Bbb C}_Q$, where $D \in \Pic^0(X)$ and
$P,Q \in X$.

For ${\cal E} \in {\bf D}(X_1 \times X_2)$, 
we define
${\cal F}_{\cal E}:{\bf D}(X_1) \to {\bf D}(X_2)$ by
\begin{equation}
{\cal F}_{\cal E}(x):={\bf R}p_{2*}({\cal E} \otimes p_1^*(x)),
x \in {\bf D}(X_1).
\end{equation}
and
$\widehat{\cal F}_{\cal E}:{\bf D}(X_2) \to {\bf D}(X_1)$ by
\begin{equation}
\widehat{\cal F}_{\cal E}(y):={\bf R}\Hom_{p_{1}}({\cal E}, p_2^*(y)),
y \in {\bf D}(X_2).
\end{equation}
We denote the $i$-th cohomology $H^i({\cal F}_{\cal E}(x))$ by
${\cal F}_{\cal E}^i(x)$.

${\cal F}_{\cal E}$ also induces isomorphisms
${\cal F}_{\cal E}:K(X_1) \to K(X_2)$,
${\cal F}_{\cal E}:{\Bbb Z} \oplus \NS(X_1) \oplus {\Bbb Z} \to
{\Bbb Z} \oplus \NS(X_2) \oplus {\Bbb Z}$
and we have a commutative diagram:
\begin{equation}
\begin{CD}
{\bf D}(X_1) @>{{\cal F}_{\cal E}}>> {\bf D}(X_2)\\
@VVV @VVV\\
K(X_1) @>{{\cal F}_{\cal E}}>> K(X_2)\\
@V{\gamma}VV @VV{\gamma}V\\
{\Bbb Z} \oplus \NS(X_1) \oplus {\Bbb Z}@>{{\cal F}_{\cal E}}>>
{\Bbb Z} \oplus \NS(X_2) \oplus {\Bbb Z}
\end{CD}
\end{equation}

For a divisor $D$,
let $T_D:{\bf D}(X) \to {\bf D}(X)$ be an equivalence of derived categories
defined by $T_D(x)=x \otimes {\cal O}_X(D)$.
This is nothing but a Fourier-Mukai transform defined by
${\cal O}_{\Delta}(D)$, where $\Delta \subset
X \times X$ is the diagonal.

\begin{lem}\label{lem:composition}
Let ${\cal F}_{{\cal E}_1}:{\bf D}(X_1) \to {\bf D}(X_2)$ and
${\cal F}_{{\cal E}_2}:{\bf D}(X_2) \to {\bf D}(X_3)$
be Fourier-Mukai transforms defined by 
${\cal E}_1 \in {\bf D}(X_1 \times X_2)$
and ${\cal E}_2 \in {\bf D}(X_2 \times X_3)$. 
Then we have
${\cal F}_{{\cal E}_2}\circ {\cal F}_{{\cal E}_1}=
{\cal F}_{{\cal G}}$, where ${\cal G}:={\cal F}_{{\cal E}_2}({\cal E}_1)$.
\end{lem}

\subsection{Twisted stability for torsion free sheaves}
Let $H$ be an ample divisor on $X$.
For $G \in K(X) \otimes{\Bbb Q}$ of $\rk G>0$, we define $G$-twisted rank, 
degree, and Euler characteristic of $x \in K(X) \otimes {\Bbb Q}$ by
\begin{equation}
\begin{split}
\rk_{G}(x)&:=\rk(G^{\vee} \otimes x)\\
\deg_{G}(x)&:=(c_1(G^{\vee} \otimes x),H)\\
\chi_{G}(x)&:=\chi(G^{\vee} \otimes x).
\end{split}
\end{equation}
For $t \in {\Bbb Q}_{>0}$, we get
\begin{equation}
\frac{\deg_G(x)}{\rk_G(x)}=\frac{\deg_{tG}(x)}{\rk_{tG}(x)},\;
\frac{\chi_G(x)}{\rk_G(x)}=\frac{\chi_{tG}(x)}{\rk_{tG}(x)}.
\end{equation}

We shall define $G$-twisted twisted stability.
\begin{defn}
Let $E$ be a torsion free sheaf on $X$.
$E$ is $G$-twisted semi-stable (resp. stable) with respect to $H$, if
\begin{equation}
\frac{\chi_G(F(nH)}{\rk_G(F)} \leq \frac{\chi_{G}(E(nH))}{\rk_{G}(E)},
n \gg 0
\end{equation}
for $0 \subsetneq F \subsetneq E$
(resp. the inequality is strict).
\end{defn} 

For a ${\Bbb Q}$-divisor $\alpha$, we define $\alpha$-twisted stability
as ${\cal O}_X(\alpha)$-twisted stability.
This is nothing but the twisted stability introduced by Matsuki and 
Wentworth \cite{M-W:1}.
It is easy to see that $G$-twisted stability is determined by
$\alpha=\det(G)/\rk G$.
Hence $G$-twisted stability is the same as the 
Matsuki-Wentworth stability.

\begin{defn}
Let ${\cal M}_{H}^{G}(\gamma)^{ss}$ be the moduli stack
of $G$-twisted semi-stable sheaves $E$ of
$\gamma(E)=\gamma$ and
${\cal M}_{H}^{G}(\gamma)^{s}$ the open substack consisting of $G$-twisted
stable sheaves.
For usual stability, {\it i.e,}
$G={\cal O}_X$,
we denote ${\cal M}_{H}^{{\cal O}_X}(\gamma)^{ss}$ by
${\cal M}_{H}(\gamma)^{ss}$. 
\end{defn} 

\begin{rem}\label{rem:beta}
Let $c_1(G)/\rk G=a H+\beta$, $a \in {\Bbb Q}, \beta \in H^{\perp}$ be
the orthogonal decomposition. Then the  
twisted semi-stability only depends on $\beta$, {\it i.e},
${\cal M}_H^{G}(\gamma)^{ss}={\cal M}_H^{\beta}(\gamma)^{ss}$.
\end{rem}

\begin{thm}\cite{M-W:1}(also see \cite{E-G:1})\label{thm:M-W}
\begin{enumerate}
\item
There is a coarse moduli scheme $\overline{M}_H^G(\gamma)$
of $S$-equivalence classes of $G$-twisted semi-stable sheaves.
\item $\overline{M}_H^G(\gamma)$ is projective.
\item For different $G,G'$, the relation between 
$\overline{M}_H^G(\gamma)$ and $\overline{M}_H^{G'}(\gamma)$
is described as Mumford-Thaddeus type flips:
\begin{equation}
\begin{matrix}
          \overline{M}_H^{G_1}(\gamma) &&&& \overline{M}_H^{G_2}(\gamma)&&&&&
    \overline{M}_H^{G_n}(\gamma)\cr
          &\searrow &&\swarrow&&\searrow&&\cdots & \swarrow &\cr
          &&\overline{M}_H^{G_{1,2}}(\gamma)&&&&\overline{M}_H^{G_{2,3}}(\gamma)&&\cr    
\end{matrix}
\end{equation}
where $G=G_1$, $G'=G_n$.
\end{enumerate}
\end{thm}

\begin{defn}
$M_H^G(\gamma)$ is the open subscheme of
$\overline{M}_H^G(\gamma)$ consisting of $G$-twisted stable sheaves and
$M_H(\gamma)^{\text{$\mu$-$s$}}$ the open subscheme 
consisting of $\mu$-stable sheaves.
Usually we denote $\overline{M}_H^{{\cal O}_X}(\gamma)$
by $\overline{M}_H(\gamma)$ and
${M}_H^{{\cal O}_X}(\gamma)$ by ${M}_H(\gamma)$.
\end{defn}
Since $\mu$-stability does not depend on $G$,
$M_H(\gamma)^{\text{$\mu$-$s$}}$ is a subscheme of $\overline{M}_H^G(\gamma)$
for all $G$.

\begin{defn}\label{defn:tw-gen}
For a pair $(H,G)$ of an ample divisor $H$ and an element
$G \in K(X) \otimes {\Bbb Q}$,
$(H,G)$ is general with respect to $v$, if
the following condition holds for every $E \in {\cal M}_H^G(v)^{ss}$: 

For $0 \subsetneq F \subsetneq E$,
\begin{equation}
\frac{\chi_G(F(nH)}{\rk_G(F)} = \frac{\chi_{G}(E(nH))}{\rk_{G}(E)},
n \gg 0
\end{equation}
implies that $v(F)/\rk F=v/\rk v$.
\end{defn}

The following is easy (cf. \cite{M-W:1}).
\begin{lem}
For an ample diviosr $H$ and a Mukai vector $v$,
there is a general $(H,G)$.
\end{lem}

\section{Fourier-Mukai transform on abelian and K3 surfaces}
\label{sect:K3}

Let $X$ be a K3 surface or an abelian surface.
Let $E$ be a coherent sheaf.
Let 
\begin{equation}
\begin{split}
v(E):=&\ch(E) \sqrt{\td_X}\\
=&\rk(E)+c_1(E)+(\chi(E)-\epsilon \rk(E))\omega_X \in H^{ev}(X,{\Bbb Z})
\end{split}
\end{equation}
be Mukai vector of $E$, where 
$\epsilon=0,1$ according as $X$ is an abelian surface or a K3 surface
and $\omega_X$ is the fundamental class of $X$.
For these surfaces, it is common to use Mukai vector
of $E$ instead of using $\gamma(E)$.
Hence we use Mukai vector in this section. 
For a Mukai vector $v$,
we define ${\cal M}(v)$, ${\cal M}_H^G(v)$, $\overline{M}_H^G(v)$,\dots
as in section 1.
By using twisted stability introduced by Matsuki and
Wentworth \cite{M-W:1},
we shall generalize \cite[sect. 8.2]{Y:7}.
In order to state our theorem (Theorem \ref{thm:00/08/17}),
we prepare some notations.

Let $v_1:=r_1+c_1+a_1 \omega_X, r_1>0, c_1 \in \NS(X)$ be a 
primitive isotropic Mukai vector on $X$.
We take a general ample divisor $H$
with respect to $v_1$.
We set $Y:=M_H(v_1)$. 
Assume that there is a universal family ${\cal E}$ on
$X \times Y$.
If $X$ is an abelian surface, then $Y$ consists of $\mu$-stable
vector bundles.
By the proof of \cite[Lem. 2.1]{Y:8}, the following lemma holds.
\begin{lem}\label{lem:class-v_1}
Assume that $H$ is general with respect to $v_1$.
\begin{enumerate}
\item
If $Y$ contains a non-locally free sheaf,
then there is an exceptional vector bundle $E_0$ and
$v_1=\rk(E_0)v(E_0)-\omega_X$.
Moreover $Y \cong X$ and a universal family is given by
\begin{equation}\label{eq:E_0}
{\cal E}:=\ker(E_0 \boxtimes E_0^{\vee} \to {\cal O}_{\Delta}).
\end{equation}
\item
If $Y$ consists of locally free sheaves, then 
they are $\mu$-stable.
\end{enumerate}
\end{lem}

We set $w_1:=v({\cal E}_{|\{x \}\times Y})=
r_1+\widehat{c}_1+\widehat{a}_1 \omega_Y$, $x \in X$. 
We consider a functor
${\cal H}_{\cal E}:
{\bf D}(X) \to {\bf D}(Y)_{op}$ defined by
\begin{equation}
 {\cal H}_{\cal E}(x):={\bf R}\Hom_{p_Y}(p_X^*(x),{\cal E}),
 x \in {\bf D}(X),
\end{equation}
where $p_X:X \times Y \to X$ (resp. $p_Y:X \times Y \to Y$) be the projection.
Then Bridgeland \cite{Br:2} showed that
${\cal H}_{\cal E}$ gives an equivalence of categories and
the inverse is given by
\begin{equation}
 \widehat{{\cal H}}_{\cal E}(y):={\bf R}\Hom_{p_X}(p_Y^*(y),{\cal E}),
 y \in {\bf D}(Y)_{op}.
\end{equation}
${\cal H}_{\cal E}$ induces an isometry $H^{ev}(X,{\Bbb Z}) \to 
H^{ev}(Y,{\Bbb Z})$.
We also denote it by ${\cal H}_{\cal E}$.

We have an isomorphism
$\NS(X) \otimes {\Bbb Q} \to v_1^{\perp} \cap \omega_X^{\perp}$
by sending $D \in \NS(X) \otimes {\Bbb Q}$ to 
$D+\frac{1}{r}(D,c_1) \omega_X \in v_1^{\perp} \cap \omega_X^{\perp}$.
Since ${\cal H}_{\cal E}$ is an isometry of Mukai lattice, 
we get an isomorphism
$v_1^{\perp} \cap \omega_X^{\perp} \to w_1^{\perp} \cap \omega_Y^{\perp}$.
Thus 
we have an isomorphism
$\delta:\NS(X) \otimes {\Bbb Q} \to \NS(Y) \otimes {\Bbb Q}$
given by
\begin{equation}
\delta(c_1(L))=c_1(p_{Y*}(\ch {\cal E} \sqrt{\td_X} 
p_X^*(c_1(L)+\frac{1}{r}(c_1(L),c_1) \omega_X)^{\vee})).
\end{equation}
For a ${\Bbb Q}$-line bundle $L \in \Pic(X) \otimes {\Bbb Q}$,
we choose a ${\Bbb Q}$-line bundle $\widehat{L}$ on $Y$ such that
$\delta(c_1(L))=c_1(\widehat{L})$.
By a result of Li \cite{Li:1} (or \cite{BBH:2}) and \cite[Lem. 7.1]{Y:7},
$\widehat{H}$ is ample, if $Y$ consists of $\mu$-stable vector bundles.
By \cite[Lem. 2.1]{Y:8}, 
$Y$ consists of $\mu$-stable vector bundles
unless ${\cal E}$ is given by \eqref{eq:E_0}.
In this case, a direct computation (or \cite{Li:1}) shows that
$\widehat{L}$ is ample.

We consider the following two conditions.
\begin{itemize}
\item[(\#1)]
$\widehat{H}$ is general with respect to $w_1$.
\item[(\#2)]
${\cal E}_{|\{x \} \times Y}$ is stable with respect to $\widehat{H}$.
\end{itemize}
\begin{rem}\label{rem:00/08/17}
If $X$ is abelian or $Y$ consists of non-locally free sheaves, then
the assumption $(\#1,2)$ holds for all general $H$.
For another example, see \cite{BBH:1}.
\end{rem}
{\bf Problem.}\cite{Y:7} Is ${\cal E}_{|\{x \} \times Y}$ always stable with respect to
$\widehat{H}$?

\vspace{1pc}

For a coherent sheaf $E$ on $X$
(resp. $F$ on $Y$), we set $\deg(E):=(c_1(E),H)$
(resp. $\deg(F):=(c_1(F),\widehat{H})$).
We consider twisted degree 
$\deg_{G_1}(E)$ and $\deg_{G_2}(F)$,
where $G_1:={\cal E}_{|X \times \{y \}}$ and
$G_2:={\cal E}_{|\{x \} \times Y}$.
Then 
\begin{lem}\cite[Lem. 8.3]{Y:7}\label{lem:degree}
$\deg_{G_1}(v)=\deg_{G_2}({\cal H}_{\cal E}(v))$. 
\end{lem}

Every Mukai vector $v$ is uniquely written as 
\begin{equation}
v=l v_1-a \omega_X+d(H+\frac{1}{r}(H,c_1)\omega_X)+(D
+\frac{1}{r}(D,c_1)\omega_X),
\end{equation}
where $l, a, d \in {\Bbb Q}$,
and $D \in \NS(X)\otimes {\Bbb Q} \cap H^{\perp}$.

It is easy to see that
$l=-\langle v,\omega_X \rangle/\rk v_1$,
$a=\langle v, v_1 \rangle/\rk v_1$ and 
$d=\deg_{G_1}(v)/r(H^2)$.
\begin{defn}
For a Mukai vector $v$,
we set
$l(v):=-\langle v,\omega_X \rangle/\rk v_1$,
$a(v):=\langle v, v_1 \rangle/\rk v_1$.
\end{defn}
Since ${\cal H}_{\cal E}(v_1)=\omega_Y$ and
$\widehat{{\cal H}}_{\cal E}(w_1)=\omega_X$,
we get
\begin{equation}\label{eq:deg-preserve}
{\cal H}_{\cal E}(l v_1-a \omega_X+(dH+D+\frac{1}{r}(dH+D,c_1)\omega_X))=
l \omega_Y-a w_1+
(d\widehat{H}+\widehat{D}+\frac{1}{r}(d\widehat{H}+\widehat{D},\widehat{c}_1)\omega_Y)
\end{equation}
where $\widehat{D} \in \NS(X) \otimes {\Bbb Q} \cap H^{\perp}$.
We can now state our theorem.

\begin{thm}\label{thm:00/08/17}
We assume the condition $(\#1,2)$ holds.
Assume that $\deg_{G_1}(v)=0$ and $l(v),a(v)>0$.
Let $\varepsilon$ be an element of $K(X) \otimes {\Bbb Q}$
such that $v(\varepsilon) \in
v_1^{\perp} \cap \omega_X^{\perp}$,
$|\langle v(\varepsilon)^2 \rangle| \ll 1$ and
$(H,c_1(\varepsilon))=0$. 
Then 
\begin{equation}
{\cal M}_H^{G_1+\varepsilon}(v)^{ss} \to 
{\cal M}_{\widehat{H}}^{G_2+\widehat{\varepsilon}}(-{\cal H}_{\cal E}(v))^{ss}.
\end{equation}
In particular, if $c_1(G_1) \in {\Bbb Q}H$, then 
$c_1(G_2) \in {\Bbb Q}\widehat{H}$ and we have an isomorphism
${\cal M}_H(v)^{ss} \to 
{\cal M}_{\widehat{H}}(-{\cal H}_{\cal E}(v))^{ss}$.
\end{thm}

The proof of Theorem \ref{thm:00/08/17} is almost the same as that in 
\cite[Thm. 8.2]{Y:7}.
Before proving Theorem \ref{thm:00/08/17},
we prepare two lemmas.

\begin{lem}\label{lem:G^2-2}
Assume that $a(v) >0$.
Then $\Hom({\cal E}_{|X \times \{y \}},E)=0$
for all $y \in Y$ and $E \in {\cal M}_H^{G_1}(v)^{ss}$.
\end{lem}

\begin{proof}
Since $H$ is general with respect to $v_1$, 
${\cal E}_{|X \times \{y \}}$ is
$G_1$-twisted stable.
Since $E$ is $G_1$-twisted semi-stable, it is sufficient to show that
$-a({\cal E}_{|X \times \{y \}})/l({\cal E}_{|X \times \{y \}})> -a(v)/l(v)$.
Since $v({\cal E}_{|X \times \{y \}})=v_1$,
we get
\begin{equation}
\frac{-a({\cal E}_{|X \times \{y \}})}
{l({\cal E}_{|X \times \{y \}})}-\frac{-a(v)}{l(v)}
=\frac{a(v)}{l(v)}>0.
\end{equation}

\end{proof} 

\begin{lem}\label{lem:G^0-2}
For a $\mu$-semi-stable sheaf $E$ of $v(E)=v$, 
there is a finite subset $S \subset Y$
such that 
\begin{equation}
\Hom(E,{\cal E}_{|X \times \{y \}})=0
\end{equation}
 for all $y \in Y \setminus S$.
\end{lem}

\begin{proof}
Considering Jordan-H\"{o}lder filtration of $E$ with respect to
$\mu$-stability,
we may assume that $E$ is $\mu$-stable.
If ${\cal E}_{|X \times \{y \}}$ is locally free, then
by Lemma \ref{lem:class-v_1}, ${\cal E}_{|X \times \{y \}}$
is $\mu$-stable, 
and hence $E^{\vee \vee} \cong {\cal E}_{|X \times \{y \}}$.
Therefore $y$ is uniquely determined by $E$.
Next we assume that ${\cal E}_{|X \times \{y \}}$ is not locally free.
Under the notation of \eqref{eq:E_0},
if $E^{\vee \vee} \ne E_0$, then clearly $\Hom(E,E_0)=0$.
Hence $\Hom(E,{\cal E}_{|X \times \{y \}})=0$ for all $y \in Y$.
If $E^{\vee \vee}=E_0$, then
$\Hom(E,{\cal E}_{|X \times \{y \}})=0$ for 
$y \in Y \setminus \Supp(E^{\vee \vee}/E)$.
\end{proof}

{\it Proof of Theorem \ref{thm:00/08/17}.}
We shall first treat the case where $\varepsilon=0$.
By the symmetry of the condition, it is sufficient to show that
$\WIT_1$ holds for $E \in M_H(v)$ 
({\it i.e,} ${\cal H}^i_{\cal E}(E)=0$, $i \ne 1$)
and
${\cal H}^1_{\cal E}(E)$ is $G_2$-twisted 
semi-stable with respect to $\widehat{L}$.
By Lemma \ref{lem:G^2-2} and
\ref{lem:G^0-2}, $\WIT_1$ holds and 
${\cal H}^1_{\cal E}(E)$ is torsion free.
We shall show that $E$ is $G_2$-twisted semi-stable.

(I) ${\cal H}^1_{\cal E}(E)$ is $\mu$-semi-stable:
Assume that ${\cal H}^1_{\cal E}(E)$ is not $\mu$-semi-stable.
Let $0 \subset F_1 \subset F_2 \subset \dots \subset 
F_s= {\cal H}^1_{\cal E}(E)$
be the Harder-Narasimhan filtration of ${\cal H}^1_{\cal E}(E)$
with respect to $\mu$-semi-stability.
We shall choose the integer $k$ which satisfies 
$\deg_{G_2}(F_i/F_{i-1})\geq 0, i \leq k$ and
$\deg_{G_2}(F_i/F_{i-1}) < 0, i > k$.
We claim that
$\widehat{\cal H}^0_{\cal E}(F_k)=0$ 
and $\widehat{\cal H}^2_{\cal E}({\cal H}_{\cal E}^1(E)/F_k)=0$.
Indeed $\deg_{G_2}(F_i/F_{i-1}) \geq 0, i \leq k$
and the $\mu$-semi-stability of $F_i/F_{i-1}$ 
imply that $\widehat{\cal H}^0_{\cal E}(F_i/F_{i-1}), i \leq k$
is of dimension 0.
Since $\widehat{\cal H}^0_{\cal E}(F_i/F_{i-1})$ is torsion free,
$\widehat{\cal H}^0_{\cal E}(F_i/F_{i-1})=0, i \leq k$.
Hence $\widehat{\cal H}^0_{\cal E}(F_k)=0$.
On the other hand, we also see that
$\widehat{\cal H}^2_{\cal E}(F_i/F_{i-1})=0, i>k$.
Hence we conclude that $\widehat{\cal H}^2_{\cal E}
({\cal H}_{\cal E}^1(E)/F_k)=0$.
So $F_k$ and ${\cal H}_{\cal E}^1(E)/F_k$ satisfy $\WIT_1$ and
we get an exact sequence 
\begin{equation}
0 \to \widehat{\cal H}^1_{\cal E}({\cal H}_{\cal E}^1(E)/F_k) \to E \to 
\widehat{\cal H}^1_{\cal E}(F_k) \to 0.
\end{equation}
By \eqref{eq:deg-preserve},
 $\deg_{G_1}(\widehat{\cal H}^1_{\cal E}(F_k))=-\deg_{G_2}(F_k)<0$.
This means that
$E$ is not $\mu$-semi-stable with respect to $L$.
Therefore ${\cal H}_{\cal E}^1(E)$ is $\mu$-semi-stable with respect to $L$.

(II) ${\cal H}^1_{\cal E}(E)$ is $G_2$-twisted semi-stable:
Assume that ${\cal H}^1_{\cal E}(E)$ is not $G_2$-twisted semi-stable.
Then there is an exact sequence
\begin{equation}
0 \to F_1 \to {\cal H}^1_{\cal E}(E) \to F_2 \to 0
\end{equation}
such that (i) $F_2$ is $G_2$-twisted stable and
(ii) $-a(F_2)/l(F_2)<-a({\cal H}^1_{\cal E}(E))/l({\cal H}^1_{\cal E}(E))
=-l(v)/a(v)$,
where $v(F_2)=l(F_2)w_1-a(F_2)\omega_Y$.
Since $ -a(F_2)/l(F_2)<-l(v)/a(v)<0$,
Lemma \ref{lem:G^2-2} and \ref{lem:G^0-2} imply that
$\widehat{{\cal H}}^0_{\cal E}(F_2)=\widehat{{\cal H}}^2_{\cal E}(F_2)=0$.
We also obtain that
$\widehat{{\cal H}}^0_{\cal E}(F_1)=0$.
Hence we have an exact sequence
\begin{equation}
 0 \to \widehat{{\cal H}}^1_{\cal E}(F_2) \to E \to 
\widehat{{\cal H}}_{\cal E}^1(F_1) \to 0.
\end{equation}
Since $\widehat{\cal H}_{\cal E}^1({\cal H}^1_{\cal E}(E))=E$,
$\widehat{\cal H}^2_{\cal E}(F_1)=0$.
Thus $\WIT_1$ also holds for $F_1$.
By (ii), we see that
\begin{equation}
\begin{split}
\frac{-a(v)}{l(v)}-\frac{-a(\widehat{\cal H}^1_{\cal E}(F_2))}
{l(\widehat{{\cal H}}_{\cal E}^1(F_2))} 
&=
\frac{-a(v)}{l(v)}+\frac{l(F_2)}{a(F_2)}\\
&=\frac{-a(v) a(F_2)+l(v) l(F_2)}{l(v)a(F_2)}
<0.
\end{split}
\end{equation}
This means that $E$ is not $G_1$-twisted semi-stable.
Therefore ${\cal H}_{\cal E}^1(E)$ is $G_2$-twisted semi-stable.

We next treat general cases.
Since $|\langle \varepsilon^2 \rangle| \ll 1$, we have an inclusion
${\cal M}_H^{G_1+\varepsilon}(v)^{ss} \subset 
{\cal M}_H^{G_1}(v)^{ss}$ and the complement consists 
of $E$ which fits in an exact sequence:
\begin{equation}\label{eq:epsilon}
0 \to E_1 \to E \to E_2 \to 0
\end{equation}
where $E_1$ is a $G_1$-twisted semi-stable sheaf such that
$v(E_1)=l_1 v_1-a_1 \omega+\delta_1$,
$\delta_1 \in v_1^{\perp} \cap \omega_X^{\perp} \cap H^{\perp}$,
$a_1/l_1=a(v)/l(v)$ and
$-\langle v(E_1),v_1+\epsilon \rangle/l_1>
-\langle v,v_1+\epsilon \rangle/l(v)$.
Then we see that 
$-\langle \delta_1,\epsilon \rangle/l_1>
-\langle \delta,\epsilon \rangle/l(v)$,
where $\delta:=v-(l(v)v_1-a(v)\omega_X) \in 
v_1^{\perp} \cap \omega_X^{\perp} \cap H^{\perp}$.
Applying ${\cal H}_{\cal E}^1$ to the exact sequence
\eqref{eq:epsilon}, we get an exact sequence
\begin{equation}
 0 \to {\cal H}^1_{\cal E}(E_2) \to {\cal H}^1_{\cal E}(E) \to 
{\cal H}_{\cal E}^1(E_1) \to 0.
\end{equation}
Since $-\langle \widehat{\delta_1},\widehat{\epsilon} \rangle/a_1>
-\langle \widehat{\delta},\widehat{\epsilon} \rangle/a(v)$,
we get that
$-\langle v({\cal H}_{\cal E}^1(E_1)), \widehat{\epsilon} \rangle/a_1<
-\langle v({\cal H}^1_{\cal E}(E)),\widehat{\epsilon} \rangle/a(v)$. 
Therefore ${\cal H}^1_{\cal E}(E)$ is not $(G_2+\widehat{\varepsilon})$-
twisted semi-stable.
\qed


\begin{prop}
Assume that $M_H(v)^{\text{$\mu$-$s$}}$ is an open dense subscheme of
$\overline{M}_H(v)$.
If $\deg_{G_1}(v)=0$, then ${\cal H}_{\cal E}$ induces a birational
map $\overline{M}_H(v) \cdots \to \overline{M}_H({\cal H}_{\cal E}(v))$
which is described as Mumford-Thaddeus type flips:
\begin{equation}
\begin{matrix}
          \overline{M}_H^{\alpha_1}(v) &&&& \overline{M}_H^{\alpha_2}(v)&&&&&
    \overline{M}_{\widehat{H}}^{\alpha_n}({\cal H}_{\cal E}(v))\cr
          &\searrow &&\swarrow&&\searrow&&\cdots & \swarrow &\cr
          &&\overline{M}_H^{\alpha_{1,2}}(v)&&&&\overline{M}_H^{\alpha_{2,3}}(v)&&\cr    
\end{matrix}
\end{equation} 
where $\alpha_i,\alpha_{i,i+1} \in \NS(X) \otimes {\Bbb Q}$ and
$\alpha_1=\alpha_n=0$. 
\end{prop}
  
\begin{proof}
By Theorem \ref{thm:M-W}, we have Mumford-Thaddeus type flips
$\overline{M}_H(v) \cdots \to \overline{M}_H^{G_1}(v)$ and 
$\overline{M}_{\widehat{H}}^{G_2}({\cal H}_{\cal E}(v)) \cdots \to 
\overline{M}_{\widehat{H}}({\cal H}_{\cal E}(v))$.
By Theorem \ref{thm:00/08/17}, $\overline{M}_H^{G_1}(v) \cong
\overline{M}_{\widehat{H}}^{G_2}({\cal H}_{\cal E}(v))$.
Therefore we get our claim.
\end{proof}

\begin{ex}
Let $X$ be a K3 surface and $H$ an ample divisor on $X$.
Assume that $H^{\perp}={\Bbb Z}D$ and $(D^2)=-2n$, $n>2$.
We set $v=2+(1-2n)\omega_X$.
Then there is a non-trivial extension
\begin{equation}
0 \to I_x(D) \to E \to {\cal O}_X(-D) \to 0
\end{equation}
where $x \in X$.
We can easily show that $E$ is a stable sheaf of $v(E)=v$.
We consider Fourier-Mukai transform defined by
${\cal E}=I_{\Delta} \otimes p_X^*{\cal O}_X(D)$.
Since $\Ext^2(E,{\cal E}_x)=\Hom({\cal E}_x,E)^{\vee} \ne 0$,
$E$ does not satisfy $\WIT_1$ with respect to ${\cal H}_{\cal E}$.
In this case, we get the following diagram
\begin{equation}
\begin{matrix}
          \overline{M}_H(v) && \cdots \to && \overline{M}_H^{D}(v) \cong
\overline{M}_H({\cal H}_{\cal E}(v))\cr
          &\searrow &&\swarrow&&\cr
          &&\overline{M}_H^{tD}(v)&&&&\cr    
\end{matrix}
\end{equation}
where $t=1/4n$.
\end{ex}

\begin{rem}
Let $(X,H)$ be a polarized K3 surface which has a divisor $\ell$ such that
\begin{equation}
(H^2)=2,\; (\ell^2)=-12,\; (H,\ell)=0
\end{equation}
and $H^0(X,{\cal O}_X(\ell+2H))=0$.
Then $Y:=M_H(2+\ell-3 \omega_X)$ is isomorphic to $X$ 
and there is a universal family ${\cal E}$ on 
$X \times Y$.
In \cite{B-M:1}, Bruzzo and Maciocia showed that
Fourier-Mukai transform
${\cal F}_{\cal E}$ gives an isomorphism
\begin{equation}
{\cal M}_H(1+(1-n)\omega_X)^{ss} \cong 
{\cal M}_{\widehat{H}}((1+2n)-n \widehat{\ell}+(1-3n)\omega_Y)^{ss}.
\end{equation}
Moreover every element $E$ of
${\cal M}_{\widehat{H}}((1+2n)-n \widehat{\ell}+(1-3n)\omega_Y)$
fits in a non-trivial extension 
\begin{equation}
0 \to E' \to E \to {\cal O}_Y \to 0
\end{equation}
where $E' \in {\cal M}_{\widehat{H}}(n(2-\widehat{\ell}-3 \omega_Y))^{ss}$.
Then we can show that $E \mapsto E^{\vee}$ induces an isomorphism
${\cal M}_{\widehat{H}}((1+2n)-n \widehat{\ell}+(1-3n)\omega_Y)^{ss} \to
{\cal M}_{\widehat{H}}^{\widehat{\ell}/2}
((1+2n)+n \widehat{\ell}+(1-3n)\omega_Y)^{ss}$.
Thus we get an isomorphism 
\begin{equation}
{\cal M}_H(1+(1-n)\omega_X)^{ss} \cong 
{\cal M}_{\widehat{H}}^{\widehat{\ell}/2}
((1+2n)+n \widehat{\ell}+(1-3n)\omega_Y)^{ss},
\end{equation}
which is nothing but the isomorphism given by ${\cal H}_{{\cal E}^{\vee}}$.
\end{rem}

\section{Irreducibility of $\overline{M}_H(v)$}\label{sect:deform}

\subsection{A special case of Theorem \ref{thm:00/08/17}}
We shall give an application of Theorem \ref{thm:00/08/17}.
Let $X$ be an abelian surface or a K3 surface such that
$\NS(X)={\Bbb Z}e \oplus {\Bbb Z}f$,
$(e^2)=(f^2)=0$ and $(e,f)=1$.
  
\begin{cor}\label{cor:ext}
Assume that (1)
$X$ is an abelian surface, $Y$ is the dual abelian surface and 
${\cal E}$ is the Poincar\'{e} line bundle
on $X \times Y$, or 
(2) $X$ is a K3 surface,
$Y=X$ and ${\cal E}$ is the ideal sheaf of the diagonal 
$\Delta \subset X \times X$
Assume that $e+kf$ is an ample divisor.
We set $D:=e-kf$.
Then ${\cal H}_{\cal E}$ induces an isomorphism of stacks
\begin{equation}
{\cal M}_{e+kf}(r+cD-a \omega_X)^{ss} \to 
{\cal M}_{\widehat{e}+k\widehat{f}}(a-c\widehat{D}-r \omega_Y)^{ss},
\end{equation}
where $r,a>0$ and $c \geq 0$.
Moreover, if $k$ is a sufficiently large integer depending
on $r$ and $\langle(r+cD-a \omega_X)^2 \rangle$, then
\begin{equation}
{\cal M}_{e+nf}(r+cD-a \omega_X)^{ss} \cong 
{\cal M}_{\widehat{e}+n\widehat{f}}(a-c\widehat{D}-r \omega_Y)^{ss}
\end{equation}
if $0<n-k \ll 1$, where $n \in {\Bbb Q}$.
\end{cor}

\begin{rem}

By Hodge index theorem,
$(D^2) < 0$.
Hence if $c>0$, then
$a> 0$, or $a=0$ and $(D^2)=-2$.
\end{rem}

\begin{proof}
By Theorem \ref{thm:00/08/17},
we get the first claim.
We next show the second claim.
We note that $(D^2)=-2k \ll 0$.
Hence $e+kf$ is a general polarization with respect to $r+cD-a \omega_X$
(cf. \cite[Lem. 2.1, Rem. 2.1]{Y:2}).
The same is true for $e+nf$, $n>k$.
Hence ${\cal M}_{e+kf}(r+cD-a \omega_X)^{ss} \cong 
{\cal M}_{e+nf}(r+cD-a \omega_X)^{ss}$.
In order to prove our claim, it suffices to show that
${\cal M}_{\widehat{e}+n\widehat{f}}(a-c\widehat{D}-r \omega_Y)^{ss}=
{\cal M}_{\widehat{e}+k\widehat{f}}(a-c\widehat{D}-r \omega_Y)^{ss}$.
We first show that
${\cal M}_{\widehat{e}+k\widehat{f}}(a-c\widehat{D}-r \omega_Y)^{ss} \subset 
{\cal M}_{\widehat{e}+n\widehat{f}}(a-c\widehat{D}-r \omega_Y)^{ss}$.
Assume that there is an exact sequence
\begin{equation}
0 \to F_1 \to {\cal H}_{\cal E}^1(E) \to F_2 \to
0 
\end{equation}
such that $F_2$ is semi-stable with respect to $\widehat{e}+n\widehat{f}$ and
\begin{enumerate}
\item
\begin{equation}
\frac{(c_1({\cal H}_{\cal E}^1(E)),\widehat{e}+n\widehat{f})}{\rk {\cal H}_{\cal E}^1(E)}>
\frac{(c_1(F_2),\widehat{e}+n\widehat{f})}{\rk F_2}
\end{equation}
or
\item
\begin{equation}
\frac{(c_1({\cal H}_{\cal E}^1(E)),\widehat{e}+n\widehat{f})}
{\rk {\cal H}_{\cal E}^1(E)}=
\frac{(c_1(F_2),\widehat{e}+n\widehat{f})}{\rk F_2},\;
\frac{\chi({\cal H}_{\cal E}^1(E))}{\rk {\cal H}_{\cal E}^1(E)}>
\frac{\chi(F_2)}{\rk F_2}.
\end{equation}
\end{enumerate}
Since $0<n-k \ll 1$, (i) or (ii) implies that
$F_1$ and $F_2$ are $\mu$-semi-stable of
$(c_1(F_1),\widehat{e}+k\widehat{f})=
(c_1(F_2),\widehat{e}+k\widehat{f})=0$ with respect to
$\widehat{e}+k\widehat{f}$.
Hence $\Hom(F_1,{\cal E}_x)=\Hom(F_2,{\cal E}_x)=0$
except finite number of points of $X$.

If (i) holds,
then $\Ext^2(F_2,{\cal E}_{|X \times \{y \}})=
\Hom({\cal E}_{|X \times \{y \}},F_2)^{\vee}=0$
for all $y \in Y$, because ${\cal E}_{|X \times \{y \}}$, $y \in Y$ 
is a stable sheaf of 
$c_1({\cal E}_{|X \times \{y \}})=0$ with respect to $e+nf$ and
$(c_1(F_2),\widehat{e}+n\widehat{f})/\rk F_2< 
(-c\widehat{D},\widehat{e}+n\widehat{f})/\rk {\cal H}_{\cal E}^1(E)=
-c(n-k)/\rk{\cal H}_{\cal E}^1(E) \leq 0$.
Therefore $F_1$ and $F_2$ satisfies $\WIT_1$ and we get an exact sequence
\begin{equation}
 0 \to \widehat{\cal H}_{\cal E}^1(F_2) \to E \to
 \widehat{\cal H}_{\cal E}^1(F_1) \to 0.
\end{equation}
By Lemma \ref{lem:degree},
\begin{equation}
 (0<) \frac{(c_1(E),e+nf)}{\rk {\cal H}_{\cal E}^1(E)}<
 \frac{(c_1(\widehat{\cal H}_{\cal E}^1(F_2)),e+nf)}{\rk F_2}.
\end{equation}
Since $E$ is semi-stable with respect to $e+kf$ and 
$(c_1(\widehat{\cal H}_{\cal E}^1(F_2)),e+kf)=0$,
$-\rk(F_2)/\rk \widehat{\cal H}_{\cal E}^1(F_2) \leq 
-\rk{\cal H}_{\cal E}^1(E)/\rk E$.
Hence we see that
\begin{equation}
 \frac{(c_1(E),e+nf)}{\rk E}<
 \frac{(c_1(\widehat{\cal H}_{\cal E}^1(F_2)),e+nf)}
 {\rk \widehat{\cal H}_{\cal E}^1(F_2)}.
\end{equation}
This implies that $E$ is not semi-stable with respect to $e+nf$.
Therefore (i) does not occur. 
If (ii) holds, then
$\frac{c_1({\cal H}_{\cal E}^1(E))}{\rk {\cal H}_{\cal E}^1(E)}=
\frac{c_1(F_2)}{\rk F_2}$.
By the proof of Theorem \ref{thm:00/08/17}, we get a contradiction.
Thus ${\cal M}_{\widehat{e}+k\widehat{f}}(a-c\widehat{D}-r \omega_Y)^{ss} 
\subset {\cal M}_{\widehat{e}+n\widehat{f}}(a-c\widehat{D}-r \omega_Y)^{ss}$.

We next show that 
${\cal M}_{\widehat{e}+n\widehat{f}}(a-c\widehat{D}-r \omega_Y)^{ss} \subset 
{\cal M}_{\widehat{e}+k\widehat{f}}(a-c\widehat{D}-r\omega_Y)^{ss}$.
Assume that there is an element 
$F \in {\cal M}_{\widehat{e}+n\widehat{f}}(a-c\widehat{D}-r \omega_Y)^{ss}
\setminus {\cal M}_{\widehat{e}+k\widehat{f}}(a-c\widehat{D}-r \omega_Y)^{ss}$.
Then we see that there is 
an exact sequence
\begin{equation}
0 \to F_1 \to F \to F_2 \to 0
\end{equation}
such that
(i) $(c_1(F),\widehat{e}+n\widehat{f})/\rk F 
\leq (c_1(F_2),\widehat{e}+n\widehat{f})/\rk F_2$,
(ii) $(c_1(F),\widehat{e}+k\widehat{f})/\rk F
=(c_1(F_2),\widehat{e}+k\widehat{f})/\rk F_2$,
(iii) $\chi(F)/\rk F>\chi(F_2)/\rk F_2$ and
(iv) $F_2$ is semi-stable with respect to $\widehat{e}+k\widehat{f}$. 
We set $\epsilon=0,1$ according as $X$ is an abelian surface or a
K3 surface as in section 2.
We note that
\begin{itemize}
\item[(a)]
${\cal E}_{|\{x\} \times Y}$, $x \in X$ is stable of 
$(c_1({\cal E}_{|\{x\} \times Y}),\chi({\cal E}_{|\{x\} \times Y}))=
(0,\epsilon)$
with respect to $\widehat{e}+n\widehat{f}$,
\item[(b)]
$F$ is stable of $(c_1(F),\widehat{e}+n\widehat{f}) \leq 0$
and $\chi(F)/\rk F=\epsilon-a/r<\epsilon$
with respect to $\widehat{e}+n\widehat{f}$.
\end{itemize}
By (a) and (b), we get
$\Ext^2(F,{\cal E}_{|\{x\} \times Y})=
\Hom({\cal E}_{|\{x\} \times Y},F)^{\vee}=0$ for all $x \in X$.
By (iii) and (iv), we see that 
$\Ext^2(F_2,{\cal E}_{|\{x\} \times Y})=0$ for all $x \in X$.
Since $F_1$ and $F_2$ are $\mu$-semi-stable sheaves of degree $0$
with respect to $\widehat{e}+k\widehat{f}$, 
$\Hom(F_1,{\cal E}_{|\{x \} \times Y})=
\Hom(F_2,{\cal E}_{|\{x \} \times Y})=0$
except finite number of points $x \in X$.
Hence $\WIT_1$ holds for $F_1$, $F_2$ and $F$ with respect to
$\widehat{{\cal H}}_{\cal E}$ and
we have an exact sequence
\begin{equation}
0 \to \widehat{{\cal H}}_{\cal E}^1(F_2) \to 
\widehat{{\cal H}}_{\cal E}^1(F) \to \widehat{{\cal H}}_{\cal E}^1(F_1)
\to 0.
\end{equation}
In the same way, we see that $E:=\widehat{{\cal H}}_{\cal E}^1(F)$ is 
$\mu$-semi-stable with respect to $e+kf$.
Since $e+kf$ is general, we get $c_1(\widehat{{\cal H}}_{\cal E}^1(F_2))/
\rk \widehat{{\cal H}}_{\cal E}^1(F_2)=c_1(\widehat{{\cal H}}_{\cal E}^1(F))/
\rk \widehat{{\cal H}}_{\cal E}^1(F)$.
On the other hand,
(i) implies that 
\begin{equation}
\frac{(c_1(E),e+nf)}{\rk F} \geq 
\frac{(c_1(\widehat{{\cal H}}_{\cal E}^1(F_2)),e+nf)}{\rk F_2}.
\end{equation}
By using (iii), we see that
\begin{equation}
\frac{(c_1(E),e+nf)}{\rk E} > 
\frac{(c_1(\widehat{{\cal H}}_{\cal E}^1(F_2)),e+nf)}{\rk 
\widehat{{\cal H}}_{\cal E}^1(F_2)},
\end{equation}
which is a contradiction.
Therefore our claim holds.
\end{proof} 

\begin{rem}
If $e+kf$ is general with respect to $r+cD-a \omega_X$ and $c>0$,
then $\Hom(E,{\cal E}_x)=0$ for all $x \in X$ and
$E \in {\cal M}_{e+kf}(r+cD-a \omega_X)^{ss}$.
Hence ${\cal H}_{\cal E}^1(E)$ is locally free.
\end{rem}

\begin{rem}
In general,
$\widehat{e}+k\widehat{f}$ is not a general polarization with respect to
$a-c\widehat{D}-r \omega_Y$. 
Indeed, let $E$ be a non-locally free $\mu$-stable sheaf of 
$v(E)=r+cD-a \omega_X$ on $X$.
Assume that $E^{\vee \vee}/E={\Bbb C}_x$, $x \in X$ and
$a >1$.
Then we get an exact sequence
\begin{equation}
0 \to {\cal H}_{\cal E}^1(E^{\vee \vee}) \to
{\cal H}_{\cal E}^1(E) \to {\cal H}_{\cal E}^2({\Bbb C}_x)
\to 0.
\end{equation}
It is easy to see that ${\cal H}_{\cal E}^2({\Bbb C}_x) \cong
{\cal E}_x$.
Hence $\widehat{e}+k\widehat{f}$ is not general with respect to
$a-c\widehat{D}-r \omega_Y$, if $c>0$.
\end{rem}

\subsection{Application to the deformation type of $M_H(v)$.}
Let $X_1$ be an abelian surface or a K3 surface.
\begin{defn}
Let $v$ be a Mukai vector of $ \rk v>0$. Then we can write it as
$v=m(v)v_p$, where $m(v) \in {\Bbb Z}$ and
$v_p$ is a primitive Mukai vector of
$\rk v_p>0$.
\end{defn}
In \cite{Y:7}, we showed that $M_H(v)$ is deformation equivalent
to a moduli space of rank 1 torsion free sheaves, if $v$ is primitive.
Here we assume that $\rk v>0$ and $H$ is general.
We shall give a slightly different proof of this result, 
that is, we shall use O'Grady's arguments \cite[sect. 2]{O:1}.
One of the benefit of O'Grady's arguments is that
we do not need to use algebraic space, which enable us to treat non-primitive
Mukai vector cases.  
For a Mukai vector $v:=l(r+c_1)+a \omega_{X_1} \in H^*(X_1,{\Bbb Z})$
such that
$r>0$, $l=\gcd(r,c_1)$ and $\gcd (l,a)=1$,
we set $b=-a+l \lambda$, $k=-(c_1^2)/2+r \lambda$, $\lambda \gg 0$
so that $e+kf$ is ample.
We consider $X$ in the above notation.
By \cite[Prop. 1.1]{Y:5} or a modification of its proof, we see that
$M_H(l(r+c_1)+a \omega_{X_1})$ is deformation equivalent to
$M_{e+nf}(l(r+(e-kf))-b \omega_X)$,
where $H$ is general with respect to $v$ and $0<n-k \ll 1$.
By Corollary \ref{cor:ext},
we have an isomorphism
\begin{equation}
M_{e+nf}(l(r+(e-kf))-b \omega_X) \to 
M_{\widehat{e}+n\widehat{f}}(b-l(\widehat{e}-k\widehat{f})-lr \omega_Y).
\end{equation}
Since $(b,l)=1$,
$M_{\widehat{e}+n\widehat{f}}(b-l(\widehat{e}-k\widehat{f})-lr \omega_Y)$ 
is deformation equivalent to
$M_{e+nf}(b+(e-k'f)-b' \omega_X)$,
where $b'=lr+\lambda', k'=l^2k+b \lambda', \lambda' \gg 0$ and
$0<n-k'\ll 1$.
Applying Corollary \ref{cor:ext} again, we get an isomorphism
\begin{equation}
M_{e+nf}(b+(e-k'f)-b' \omega_X) \to 
M_{\widehat{e}+n\widehat{f}}(b'-(\widehat{e}-k'\widehat{f})-b \omega_Y).
\end{equation}
If $\lambda'$ is sufficiently large, then
$M_{\widehat{e}+n\widehat{f}}(b'-(\widehat{e}-k'\widehat{f})-b \omega_Y)$ is 
deformation equivalent to
$M_{\widehat{e}+\widehat{n}f}(b'-(\widehat{e}-k''\widehat{f})-\omega_Y)$ 
and $e+k''f$ is ample, where
$k''=lr(1-b)+l^2k+ \lambda'$ and $0<n-k \ll 1$.
Since $k'' \gg 0$,
Corollary \ref{cor:ext} implies that
$M_{{e}+n{f}}(1+({e}-k''{f})-b'' \omega_X)$ 
is isomorphic to
$M_{\widehat{e}+\widehat{n}f}(b'-(\widehat{e}-k''\widehat{f})-\omega_Y)$.
Therefore $M_H(v)$ is deformation equivalent to 
the moduli space of rank 1 torsion free sheaves.

We shall next treat non-primitive Mukai vector.

\begin{lem}\label{lem:pss}
Let $v$ be a Mukai vector of $\rk v>0$.
Let $H$ be a general ample divisor with rerspect to $v$.
We set
\begin{equation}
 {\cal M}_H(v)^{pss}:=\{E \in {\cal M}_H(v)^{ss}| 
 \text{ $E$ is properly semi-stable }\}.
\end{equation}
Then $\dim {\cal M}_H(v)^{pss} \leq \langle v^2 \rangle$.
Moreover inequality is strict, unless $m(v)=2$ and $\langle v^2 \rangle=8$.
\end{lem}
For the proof, see \cite[Lem. 1.7]{Y:8}.
  
\begin{prop}\label{prop:normal}
Under the same assumptions,
${\cal M}_H(v)^{ss}$ is a locally complete intersection stack
which contains ${\cal M}_H(v)^{s}$ as an open dense substack 
and the singular locus is at least of codimension 2.
In particular ${\cal M}_H(v)^{ss}$ is normal.
\end{prop}

\begin{proof}
By deformation theory and Lemma \ref{lem:pss},
we see that ${\cal M}_H(v)^{ss}$ is a locally complete intersection.
If $m(v) \ne 2$ or $\langle v^2 \rangle>8$, then  
$\dim {\cal M}_H(v)^{pss} \leq \langle v^2 \rangle-1$.
Therefore the singular locus is at least of codimension 2.
If $m(v)=2$ and $\langle v^2 \rangle=8$, then
a general member of ${\cal M}_H(v)^{pss}$ fits in a non-trivial
extension
\begin{equation}
 0 \to E_1 \to E \to E_2 \to 0
\end{equation}
where $E_1, E_2 \in {\cal M}_H(v/2)^{ss}$ and $E_1 \ne E_2$.
Then $E$ is simple, which implies that
${\cal M}_H(v)^{ss}$ is smooth at $E$.
Therefore the singular locus is at least of codimension 2.
For the last claim, we use Serre's criterion.  
\end{proof}

\begin{defn}
Let $Y_1, Y_2$ be normal schemes.
Then $Y_1 \sim Y_2$,
if there is a proper and flat morphism ${\cal Y} \to T$
over a smooth connected curve $T$ such that every fiber is normal and
$Y_i={\cal Y}_{t_i}$ for some $t_1,t_2 \in T$.
Deformation equivalence is an equivalence relation generated by 
$\sim$.
\end{defn}

By the following lemma, the number of irreducible components is 
an invariant of this equivalence relation.

\begin{lem}
Let $T$ be a smooth curve and $Y \to T$ a flat and proper morphism.
Assume that every fiber is normal. Then
the number of irreducible components of $Y_t$, $t \in T$ is constant.
\end{lem}

\begin{proof}
Since $Y_t$ is normal, every connected component is an integral scheme.
Hence the number of irreducible components of $Y_t$ is 
$h^0(Y_t, {\cal O}_{Y_t})$.
By the upper-semicontinuity of $h^0(Y_t, {\cal O}_{Y_t})$,
the number of irreducible components of $Y_t$ is upper semi-continuous.
On the other hand, by Zariski's connectivity theorem,
the number of connected components of $Y_t$ is lower semi-continuous.
Therefore we get our lemma.
\end{proof}

By the same proof, we can show the following. 
\begin{prop}
$\overline{M}_H(v)$
is deformation equivalent to
$\overline{M}_H(m(v)(1-n \omega_X))$, where
$n=\langle v_p^2 \rangle/2$.
In particular the number of irreducible components of $\overline{M}_H(v)$
is determined by $m(v)$.
\end{prop}
\begin{proof}
Let $T$ be a smooth curve over ${\Bbb C}$ and
$\varphi:({\cal X},{\cal L}) \to T$ be a family of polarized
abelian or K3 surfaces. 
For a family of Mukai vectors 
$v \in R^* \varphi_* {\Bbb Z}=\cup_{t \in T}H^*({\cal X}_t,{\Bbb Z})$,
let $\psi:\overline{{\frak M}}_{\cal L}(v) \to T$ be the relative moduli space
of semi-stable sheaves on ${\cal X}_t$, $t \in T$ of
Mukai vector $v_t$
and ${\frak M}_{\cal L}(v)$ the open subscheme
consisting of stable sheaves.
Since $T$ is defined over a field of characteristic 0,
$\overline{{\frak M}}_{\cal L}(v)_t=\overline{M}_{{\cal L}_t}(v_t)$ 
for $t \in T$, where $\overline{M}_{{\cal L}_t}(v_t)$ is the moduli space of
semi-stable sheaves on ${\cal X}_t$ (cf. \cite[Thm. 1.1]{MFK:1}).
Since $\psi_{|{\frak M}_{\cal L}(v)}:
{\frak M}_{\cal L}(v) \to T$ is smooth \cite{Mu:3}, it is flat.
Assume that ${\cal L}_t$ is general with respect to $v_t$ for all 
$t \in T$.
By Proposition \ref{prop:normal},
 ${\frak M}_{\cal L}(v)$ is a dense subscheme of 
$\overline{{\frak M}}_{\cal L}(v)$.
Since $T$ is a smooth curve, $\psi$ is also flat.
Therefore $\psi:\overline{{\frak M}}_{\cal L}(v) \to T$ is a proper and flat
morphism. 
Then our claim follows from the same argument as
in $m(v)=1$ case.  
\end{proof}

\subsection{Irreducibility of $\overline{M}_H(v)$}

We shall show that $\overline{M}_H(v)$ is irreducible.
We may assume that $X$ has an elliptic fibration $\pi:X \to C$.
We also assume that there is a section $\sigma$ of $\pi$
and $\NS(X)={\Bbb Z}\sigma \oplus {\Bbb Z}f$, where $f$ is a fiber of $\pi$. 
By \cite[Thm. 3.15]{Y:7} (cf. Theorem \ref{thm:FM}), we have an isomorphism
$\overline{M}_{\sigma+kf}(r(1-n \omega_X)) \cong 
\overline{M}_{\sigma+kf}(w)$, where
$w=r((\sigma+(n+\epsilon)f)+\omega_X)$.
Hence it is sufficient to show that 
$\overline{M}_{\sigma+kf}(w)$ is irreducible.
From now on, we assume that $r \geq 2$.
By Proposition \ref{prop:normal}, we shall show that $M_{\sigma+kf}(w)$
is irreducible.

\begin{defn}
For a purely 1-dimensional sheaf $L$ on $X$,
$\Div(L)$ is the divisor on $X$ which is defined by the fitting ideal
of $L$.
\end{defn}

We set $\xi=r(\sigma+(n+\epsilon) f)$.
Let $\Hilb_X^{\xi}$ be the Hilbert scheme of curves $C$ on $X$
such that $c_1({\cal O}_X(C))=\xi$. 
There is a natural map $j:\overline{M}_{\sigma+kf}(w) \to \Hilb_X^{\xi}$
sending $L \in \overline{M}_{\sigma+kf}(w)$ to $\Div(L)$.
We want to estimate the dimension of closed subsets of 
$M_{\sigma+kf}(w)$:
\begin{equation}
\begin{split}
N_1:&=\{L \in M_{\sigma+kf}(w)|\text{ $\Div(L)$ is not irreducible}\},\\
N_2:&=\{L \in M_{\sigma+kf}(w)|\text{ $\Div(L)$ is not reduced}\}.
\end{split}
\end{equation}
{\it Estimate of $\dim N_1$.}
We prepare some lemmas. 

\begin{lem}\label{lem:intersection}
Let $C_i$, $i=1,2$ be irreducible curves of genus $g(C_i) \geq 2$.
Then $(C_i,C_j) \geq 2$.
\end{lem}
\begin{proof}
If $C_1-C_2$ or $C_2-C_1$ is effective, then
$(C_1,C_2) \geq (C_2^2) \geq 2$ or $(C_1,C_2) \geq (C_1^2) \geq 2$.
If $C_1-C_2$ and $C_2-C_1$ are not effective, then
$0 \geq \chi({\cal O}_X(C_1-C_2))\geq (C_1-C_2)^2/2$.
Hence we see that $(C_1,C_2) \geq 2$.
\end{proof}

\begin{defn}
For a Mukai vector $v \in H^{ev}(X,{\Bbb Z})$,
${\cal M}(v)$ is the stack of coherent sheaves $E$ of $v(E)=v$.
\end{defn}

\begin{lem}\label{lem:g>1}
Let $E$ be a purely 1-dimensional sheaf such that $\Supp E$ consists of 
genus $g \geq 2$ curves.
Then 
 $\dim {\cal M}(v(E))=\langle v(E)^2 \rangle+1$ at $E$, if $H$ is general.
\end{lem}

\begin{proof}
We set $v:=v(E)$.
Let ${\cal M}$ be an irreducible component of ${\cal M}(v)$ containing
$E$ and let $E'$ be a general point of ${\cal M}$.
We consider the Harder-Narasimhan filtration of $E'$:
\begin{equation}\label{eq:filtration}
 0 \subset F_1 \subset F_2 \subset \dots \subset F_s=E'.
\end{equation}
We set $v_i=v(F_i/F_{i-1})$.
By our assumption, we may assume that
$\Supp F_i/F_{i-1}$ consist of curves of genus greater than 1.
Hence $\langle v_i^2 \rangle >0$.
Moreover by Lemma \ref{lem:intersection}, $\langle v_i,v_j \rangle \geq 2$.
Let ${\cal F}^0(v_1,v_2,\dots,v_s)$ be the stack of filtrations
\eqref{eq:filtration} such that
$\Hom(F_i/F_{i-1},F_j/F_{j-1})=0$ for $i<j$.
By \cite[Lem. 5.2]{Y:8}, 
\begin{equation}
 \dim {\cal F}^0(v_1,v_2,\dots,v_s)=
 \sum_{i=1}^s \dim {\cal M}_{\sigma+kf}(v_i)^{ss}+
 \sum_{i<j}\langle v_i,v_j \rangle.
\end{equation}
Since $\langle v_i^2 \rangle>0$,
we get $\dim {\cal M}_{\sigma+kf}(v_i)^{ss}=\langle v_i^2 \rangle+1$.
Hence we see that
\begin{equation}\label{eq:v1,,,vs}
\begin{split}
\dim {\cal F}^0(v_1,v_2,\dots,v_s)&=
\langle v^2 \rangle+1-(
\sum_{i=1}^s (\langle v_i \rangle+1)+\sum_{i<j}\langle v_i,v_j \rangle)\\
&=\sum_{i<j}\langle v_i,v_j \rangle-(s-1) >0.
\end{split} 
\end{equation}
Since $\dim {\cal M}(v) \geq \langle v^2 \rangle+1$, we get our claim.
\end{proof}

\begin{lem}\label{lem:g=0}
Assume that $X$ is a K3 surface.
Let $E$ be a purely 1-dimensional sheaf of $\Div(E)=r \sigma$.
Then
${\cal M}(v(E))=-r^2$ at $E$.
\end{lem}

\begin{proof}
We set $v:=v(E)$.
Let ${\cal M}$ be an irreducible component of ${\cal M}(v)$ containing
$E$ and let $E'$ be a general point of ${\cal M}$.
We consider the Harder-Narasimhan filtration of $E'$:
\begin{equation}\label{eq:filtration2}
 0 \subset F_1 \subset F_2 \subset \dots \subset F_s=E'.
\end{equation}
We set $v_i=v(F_i/F_{i-1})$.
Then $v_i=r_i \sigma+a_i \omega_X$.
It is easy to see that $a_i$ is divisible by $r_i$ and
${\cal M}_{\sigma+kf}(v_i)^{ss}=
\{{\cal O}_{\sigma}(a_i/r_i-1)^{\oplus r_i} \}$.
Then $\dim {\cal M}_{\sigma+kf}(v_i)^{ss}=-r_i^2$.
As in Lemma \ref{lem:g>1}, let
 ${\cal F}^0(v_1,v_2,\dots,v_s)$ be the stack of filtrations
\eqref{eq:filtration2} such that
$\Hom(F_i/F_{i-1},F_j/F_{j-1})=0$ for $i<j$.
By \cite[Lem. 5.2]{Y:8}, 
\begin{equation}
\begin{split}
 \dim {\cal F}^0(v_1,v_2,\dots,v_s)&=
\sum_{i=1}^s \dim {\cal M}_{\sigma+kf}(v_i)^{ss}+
 \sum_{i<j}\langle v_i,v_j \rangle\\
 &=-\sum_{i=1}^s r_i^2-\sum_{i<j} 2r_ir_j\\
 &=-r^2.
\end{split}
\end{equation}
Therefore we get our claim.
\end{proof}

\begin{lem} \label{lem:g=1}
Let $E$ be a purely 1-dimensional sheaf on $X$
such that $v(E)=rf+a \omega_X$, or
$v(E)=r\sigma+a \omega_X$.
Assume that $\langle v(E)^2 \rangle=0$.
Then $\dim {\cal M}(v(E))= r$ at $E$.
\end{lem}

\begin{proof}
We set $v:=v(E)$.
Let ${\cal M}$ be an irreducible component of ${\cal M}(v)$ containing
$E$ and let $E'$ be a general point of ${\cal M}$.
We consider the Harder-Narasimhan filtration of $E'$:
\begin{equation}\label{eq:filtration3}
 0 \subset F_1 \subset F_2 \subset \dots \subset F_s=E'.
\end{equation}
We set $v_i=v(F_i/F_{i-1})$.
Then we see that $\langle v_i^2 \rangle=0$.
As in the proof of Lemma 1.8
in \cite{Y:8},
we see that $\dim {\cal M}_{\sigma+kf}(v_i)^{ss}=r_i$.
By using \cite[Lem. 5.2]{Y:8} again, we see that 
\begin{equation}
\begin{split}
 \dim {\cal F}^0(v_1,v_2,\dots,v_s)&=
\sum_{i=1}^s \dim {\cal M}_{\sigma+kf}(v_i)^{ss}+
 \sum_{i<j}\langle v_i,v_j \rangle\\
 &=\sum_{i=1}^s r_i=r.
\end{split}
\end{equation}
Therefore we get our claim.
\end{proof}

For $N_1$, we get the following.
\begin{prop}\label{lem:irred}
$\dim N_1<\dim M_{\sigma+kf}(w)$.
\end{prop}
\begin{proof}

Assume that $\Supp(L)$ is not
irreducible.
Then there is a filtration
\begin{equation}
 0 \subset F_1 \subset F_2 \subset F_3=L
\end{equation}
such that (i) $\Div (F_1)=r_1 \sigma$
(ii) $F_2/F_1$ is a pure dimension 1 sheaf of
$\Div (F_2/F_1)=r_2 f$ and (iii)
$F_3/F_2$ is a pure dimension 1 sheaf 
of $\Div(F_3/F_2)=C_3$,
where $C_3$ consists of curves of genus greater than 1. 
We set $v_i:=v(F_i/F_{i-1})$, $i=1,2,3$.
By Lemma \ref{lem:g>1}, we may assume that $F_2 \ne 0$.
We first note that $\Ext^2(F_i/F_{i-1},F_j/F_{j-1})=
\Hom(F_j/F_{j-1},F_i/F_{i-1})^{\vee}=0$ for $i \ne j$.

(I) We first treat the case where $X$ is a K3 surface.
Assume that $F_2 \ne F_3$.
By Lemmas \ref{lem:g>1}, \ref{lem:g=0}, \ref{lem:g=1}, we see that
\begin{equation}
\begin{split}
\dim {\cal F}^0(v_1,v_2,v_3)&=\sum_i \dim {\cal M}_{\sigma+kf}(v_i)^{ss}
+\sum_{i<j} \langle v_i,v_j \rangle\\
&=\langle w^2 \rangle+1-
(-r_1^2-r_2+(r_1\sigma,r_2f+C_3)+(r_2f,C_3)).
\end{split}
\end{equation}
By our assumption,
$(c_1(w),\sigma) \geq 0$.
Hence $(r_1\sigma,r_2f+C_3) \geq 2r_1^2$.
By our assumption, $(f,C_3)>0$.
Therefore $-r_1^2-r_2+(r_1\sigma,r_2f+C_3)+(r_2f,C_3) \geq r_1^2>0$.
We next assume that $F_2=F_3$.
Then we see that 
\begin{equation}
\dim {\cal F}^0(v_1,v_2)
=\langle w^2 \rangle+1-
(-r_1^2-r_2+1+r_1r_2).
\end{equation}
Since $(c_1(w),\sigma) \geq 0$, $r_2 \geq 2r_1$.
Then $-r_1^2-r_2+1+r_1r_2 \geq r_1(r_1-2)+1 \geq 1$,
because $c_1(w)$ is not primitive.
Therefore we get our claim.

(II) We next treat the case where $X$ is an abelian surface.
Assume that $F_2 \ne F_3$. Then
\begin{equation}
\dim {\cal F}^0(v_1,v_2,v_3)=\langle w^2 \rangle+1-
(-r_1-r_2+(r_1\sigma,r_2f+C_3)+(r_2f,C_3)).
\end{equation}
Since $C_3$ consists of curves of genus greater that 1,
$(\sigma,C_3)>1$ and $(f,C_3)>1$. 
Then $(-r_1-r_2+(r_1\sigma,r_2f+C_3)+(r_2f,C_3))>0$.
If $F_2=F_3$, then 
\begin{equation}
\dim {\cal F}^0(v_1,v_2)=\langle w^2 \rangle+1-(-r_1-r_2+1+r_1r_2).
\end{equation}
Since $c_1(w)$ is not primitive,
$(r_1-1)(r_2-1)>0$.
Therefore we get our claim.  
\end{proof}
{\it Estimate of $\dim N_2$.}
For an integer $\lambda \geq 2$, 
let $\Hilb_X^{\xi}(\lambda)$ be the locally closed subset of 
$\Hilb_X^{\xi}$ consisting of $\lambda C$, where
$C$ is an integral curve.
By Proposition \ref{lem:irred}, it is sufficient to 
estimate the dimension of 
$\Hilb_X^{\xi}(\lambda)$. 
Let $C$ be an integral curve. 
For $D=\lambda C$ and $w=D+r \omega_X$, we set
\begin{equation}
 M_{\sigma+kf}(w,D):=\{L \in M_{\sigma+kf}(w)|\Div(L)=D \}.
\end{equation}
We fix a point $x \in C$.
We also set
\begin{equation}
P_{\sigma+kf}(w,w+\omega_X,D):=\left\{L \subset L'\left|
\begin{split}
&L' \in M_{\sigma+kf}(w+\omega_X),
\Div(L')=D\\
&L \in M_{\sigma+kf}(w), L'/L \cong {\Bbb C}_x
\end{split}
\right. \right\}.
\end{equation}
Let $\pi_{w}$ and $\pi_{w+\omega_X}$ be natural projections sending
$L \subset L'$ to $L$ and $L'$ respectively:
\begin{equation}
\begin{matrix}
          && P_{\sigma+kf}(w,w+\omega_X,D) && \cr
          &\llap{$\scriptstyle{\pi_{w}}$}\swarrow&&\searrow
          \rlap{$\scriptstyle{\pi_{w+\omega_X}}$}&&\cr
          M_{\sigma+kf}(w,D)&&&&M_{\sigma+kf}(w+\omega_X,D)\cr    
\end{matrix}
\end{equation}
For $L \in M_{\sigma+kf}(w,D)$, 
$-\dim \Ext^1({\Bbb C}_x,L)+\dim \Ext^2({\Bbb C}_x,L)=
\chi({\Bbb C}_x,L)=0$.
Combining Serre duality, we see that
$\dim \Ext^1({\Bbb C}_x,L)=\dim \Ext^2({\Bbb C}_x,L)=
\dim \Hom(L,{\Bbb C}_x)$.
By the following lemma, $\dim \pi_w^{-1}(L) \leq \lambda-1$
and $\dim \pi_{w+\omega_X}^{-1}(L') \leq \lambda-1$.
\begin{lem}
Let $x$ be a smooth point of $C$.
Let $L$ be a purely 1-dimensional sheaf such that 
$\Div(L)=\lambda C$.
Then $\dim L \otimes {\Bbb C}_x \leq \lambda$.
\end{lem}

\begin{proof}
Let $C'$ be a germ of a curve intersecting $C$ at $x$ transversely.
Let ${\cal O}_{X,x}$ be the stalk of ${\cal O}_X$ at $x$.
We take a free resolution of $L \otimes {\cal O}_{X,x}$:
\begin{equation}
0 \to {\cal O}_{X,x}^{\oplus n} \overset{A}{\to} {\cal O}_{X,x}^{\oplus n}
\to L \otimes {\cal O}_{X,x} \to 0.
\end{equation}
Then the local equation of $\Div(L)$ at $x$ is given by $\det(A)$. 
By restricting the sequence to $C'$, we get a free resolution of 
$L \otimes {\cal O}_{C',x}$.
Then $\dim (L \otimes {\cal O}_{C',x})$ is given by 
the local intersection number
$(\Div(D),C')_x=\lambda$. Therefore we get our claim.
\end{proof}

\begin{lem}\label{lem:lag}
$\dim M_{\sigma+kf}(w+\omega_X,D)=(D^2)/2+1$.
\end{lem}
\begin{proof}
By \cite[Thm. 3.15]{Y:7},
$M_{\sigma+kf}(w+\omega_X)$ is isomorphic to
$M_{\sigma+kf}(r+f-rn \omega_X)$.
Since $r+f-rn \omega_X$ is primitive,
 \cite[Thm. 0.1 and 8.1]{Y:7} implies that it is irreducible.
For a smooth curve $C \in \Hilb_X^{\xi}$,
the fiber of 
$M_{\sigma+kf}(w+\omega_X) \to \Hilb_X^{\xi}$ is
$\Pic^{r+1}(C)$.
It is easy to see that $\Pic^{r+1}(C)$ is a 
Lagrangian subscheme of $M_{\sigma+kf}(w+\omega_X)$.
By Matsushita \cite{Mat:1},
every fiber is of dimension $(\xi^2)/2+1$.
\end{proof}

\begin{lem}\label{lem:L'}
Let $L$ be a stable sheaf of $v(L)=w$ and
$\Div(L)=D$, and let $L'$ be a coherent
sheaf which fits in a non-trivial extension
\begin{equation}
0 \to L \to L' \to {\Bbb C}_x \to 0
\end{equation}
where $x \in D$.
Then $L'$ is stable.
\end{lem}

\begin{proof}
Assume that $L'$ is not of pure dimension 1 and 
let $T$ be the 0-dimensional subsheaf of $L'$.
Then $T \to L' \to {\Bbb C}_x$ must be injective.
Hence it is isomorphic, which implies that the exact sequence
split.
Therefore $L'$ is of pure dimension 1.
If $L'$ is not stable, then there is a subsheaf $L_1$ of $L'$
such that $\frac{\chi(L_1)}{\lambda_1}>\frac{\chi(L')}{\lambda}$,
where $\Div(L_1)=\lambda_1 C$.
Hence $0<\chi(L_1)\lambda-\chi(L')\lambda_1=
(\chi(L_1)-1)\lambda-\chi(L)\lambda_1+\lambda-\lambda_1$.
Since $\chi(L)=r$ is divisible by $\lambda$ and
$\lambda-\lambda_1<\lambda$, we get
$(\chi(L_1)-1)\lambda-\chi(L)\lambda_1 \geq 0$.
Thus $\frac{\chi(L_1 \cap L)}{\lambda_1}\geq \frac{\chi(L)}{\lambda}$,
which implies that $L$ is not stable.
Therefore $L'$ must be stable.
\end{proof}

\begin{cor}
$\dim M_{\sigma+kf}(w,D) \leq \dim P_{\sigma+kf}(w,w+\omega_X,D) \leq 
\dim M_{\sigma+kf}(w+\omega_X,D)+(\lambda-1)$.
\end{cor}

\begin{proof}
By Lemma \ref{lem:L'}, $\pi_w$ is surjective.
Hence we get our claim from Lemma \ref{lem:lag}.
\end{proof}  
Since $\lambda \geq 2$, we get that $2\lambda^2-(\lambda^2+\lambda+1)>0$.
Then 
\begin{equation}
\begin{split}
 \dim j^{-1}(\Hilb_X^{\xi}(\lambda)) &\leq 
 (C^2)/2+1+\lambda^2(C^2)/2+1+(\lambda-1)\\
 &< (\lambda^2+\lambda+1)(C^2)/2+2\\
 &< \lambda^2(C^2)+2=\dim M_{\sigma+kf}(w).
\end{split}
\end{equation}
Combining Proposition \ref{lem:irred}, we get the following. 
\begin{prop}\label{lem:integral}
We set 
\begin{equation}
\begin{split}
 M_{\sigma+kf}(w)_0:&=\{L \in M_{\sigma+kf}(w)|
 \text{ $\Div(L)$ is an integral curve }\}\\
 &=M_{\sigma+kf}(w) \setminus(N_1 \cup N_2).
\end{split}
\end{equation}
Then $M_{\sigma+kf}(w)_0$ is an open dense subscheme of
$M_{\sigma+kf}(w)$.
\end{prop}

\begin{prop}
$M_{\sigma+kf}(w)$ is irreducible.
\end{prop}

\begin{proof}
By Proposition \ref{lem:integral},
it is sufficient to show that $M_{\sigma+kf}(w)_0$ is irreducible.
Let $C$ be an integral curve. 
Then $j^{-1}(C)$ is the compactified Jacobian of $C$.
By \cite{A-I-K:1}, the compactified Jacobian of $C$
is irreducible.
Therefore $M_{\sigma+kf}(w)_0$ is irreducible.
\end{proof}

Combining all together, we get the following theorem.
\begin{thm}
Let $X$ be an abelian surface or a K3 surface and 
let $v$ be a Mukai vector of
$\rk v>0$ and $\langle v^2 \rangle>0$.
Let $H$ be a general ample divisor with respect to $v$.
Then ${\cal M}_H(v)^{ss}$ is a normal and irreducible stack.
In particular, $\overline{M}_H(v)$ is a normal variety. 
\end{thm}

\begin{rem}
In \cite{O:2}, \cite{O:3},
O'Grady studied the case where $m(v)=2$.
In particular, he constructed symplectic desinguralization 
of $\overline{M}_H(v)$, if $\langle v^2 \rangle=8$.
\end{rem}

\section{Fourier-Mukai transform on Enriques surfaces}\label{sect:enriques}

In this section, we consider Fourier-Mukai transform 
on Enriques surface $X$.
By using Fourier-Mukai transform, we shall compute Hodge
polynomial of some moduli spaces of sheaves.

In our case, Mukai vector $v(x)$ of $x \in K(X)$
is defined as an element of $H^*(X,{\Bbb Q})$:
\begin{equation}
\begin{split}
v(x):=&\ch(x)\sqrt{\td_X}\\
=&\rk(x)+c_1(x)+(\frac{\rk(x)}{2}\omega_X+\ch_2(x))
\in H^*(X,{\Bbb Q}).
\end{split}
\end{equation}
We also introduce Mukai's pairing on $H^*(X,{\Bbb Q})$ by
$\langle x,y \rangle:=-\int_X x^{\vee} \wedge y$.
Then we have an isomorphism of lattices:
\begin{equation}
(v(K(X)),\langle \;\;,\;\;\rangle) \cong
\begin{pmatrix}
1&0\\
0&-1
\end{pmatrix}
\oplus 
\begin{pmatrix}
0&1\\
1&0
\end{pmatrix}
\oplus E_8(-1).
\end{equation}
\begin{defn}
We call an element of $v(K(X))$ by Mukai vector.
A Mukai vector $v$ is primitive, if
$v$ is primitive as an element of $v(K(X))$.
\end{defn}

The following was essentially proved in \cite[sect. 3.2]{Y:2}.
\begin{prop}\label{prop:indep}
Let $X$ be a surface such that $K_X$ is
numerically trivial.
Let $(H,\alpha)$ be a pair of ample divisor $H$
and a ${\Bbb Q}$-divisor $\alpha$.
Then
$e({\cal M}_H^{\alpha}(v)^{ss})$ does not depend on the choice of
$H$ and $\alpha$, if $(H,{\cal O}(\alpha))$ 
is general with respect to $v$ (cf. Defn. \ref{defn:tw-gen}).
\end{prop}

\begin{prop}
Let $v$ be a Mukai vector such that $\rk(v)$ is odd.
Then ${\cal M}_H(v)^s$ is smooth of 
$\dim {\cal M}_H(v)^s=\langle v^2 \rangle+1$.
\end{prop}

\begin{proof}
For $E \in {\cal M}_H(v)^s$, we get 
$\det(E(K_X)) \not \cong \det(E)$.
If there is a non-zero homomorphism
$E \to E(K_X)$, then the stability implies that
it is an isomorphism.
Hence $\Ext^2(E,E)=\Hom(E,E(K_X))^{\vee}=0$.
Since $-\chi(E,E)=-\langle v(E),v(E) \rangle$,
${\cal M}_H(v)^s$ is smooth of 
$\dim {\cal M}_H(v)^s=\langle v^2 \rangle+1$.
\end{proof}

For a Mukai vector $v$,
let $L_1, L_2=L_1(K_X) \in \Pic(X)$ be line bundles on $X$ such that
$c_1(L_1)(=c_1(L_2))=c_1(v)$.
Then we have a decomposition
\begin{equation}
{\cal M}_H(v)^{ss}={\cal M}_H(v,L_1)^{ss} \coprod
 {\cal M}_H(v,L_2)^{ss}
\end{equation}
where ${\cal M}_H(v,L_i)^{ss}$, $i=1,2$, 
is the substack of ${\cal M}_H(v)^{ss}$
consisting of $E$ such that $\det(E)=L_i$.
We also have a decomposition $M_H(v)=M_H(v,L_1) \coprod M_H(v,L_2)$,
where $M_H(v,L_i)$ is the subscheme of $M_H(v)$
consisting of $E$ such that $\det(E)=L_i$.

We consider Fourier-Mukai transform associated to $(-1)$-reflection.
Let $v_0:=r+c_1-(s/2) \omega_X$ be a Mukai vector such that
$\rk(v_0)>0$ and $\langle v_0^2 \rangle=(c_1^2)+rs=-1$.
Since $(c_1^2)$ is even,
$r$ and $s$ are odd.
Let $H$ be a general ample divisor with respect to $v_0$.
Assume that there is a stable vector bundle $E_0$
with respect to $H$ such that $v(E_0)=v_0$
(cf. Corollary \ref{cor:excep}).
Then we see that
\begin{equation}
 \begin{cases}
  \Hom(E_0,E_0)={\Bbb C}\\
  \Ext^1(E_0,E_0)=0\\
  \Ext^2(E_0,E_0)=0.
 \end{cases}
\end{equation}
Let 
\begin{equation}
\begin{split}
ev_1:&E_0^{\vee} \boxtimes E_0 \to {\cal O}_{\Delta},\\
ev_2:&E_0(K_X)^{\vee} \boxtimes E_0(K_X) \to {\cal O}_{\Delta}
\end{split}
\end{equation}
be evaluation maps.
We define a sheaf ${\cal E}$ on $X \times X$ by an exact sequence
\begin{equation}
0 \to {\cal E} \to E_0^{\vee} \boxtimes E_0 \oplus
E_0(K_X)^{\vee} \boxtimes E_0(K_X)
\overset{(ev_1,ev_2)}{\to} {\cal O}_{\Delta} \to 0. 
\end{equation}
Then ${\cal E}_{|\{x \} \times X}$
(resp. ${\cal E}_{|X \times \{x \}}$) is a stable sheaf
of $v({\cal E}_{|\{x \} \times X})=2\rk(E_0)v(E_0)-\omega_X$
(resp. $v({\cal E}_{|X \times \{x \}})=2\rk(E_0)v(E_0)^{\vee}-\omega_X$).
Thus ${\cal E}$ is a flat family of stable sheaves
of $v({\cal E}_{|\{x \} \times X})=2\rk(E_0)v(E_0)-\omega_X$.
By the construction of ${\cal E}$,
${\cal E}_{|\{x \} \times X}(K_X) \cong {\cal E}_{|\{x \} \times X}$,
which implies that
\begin{equation}
 \Ext^2({\cal E}_{|\{x \} \times X},{\cal E}_{|\{x \} \times X})=
 \Hom({\cal E}_{|\{x \} \times X},{\cal E}_{|\{x \} \times X}(K_X))^{\vee}
 \cong {\Bbb C}.
\end{equation}
Since $\langle v(E_0)^2 \rangle=-1$, we see that
$\langle v({\cal E}_{|\{x \} \times X})^2 \rangle=0$.
Hence the Zariski tangent space is 2-dimensional:
\begin{equation}
\Ext^1({\cal E}_{|\{x \} \times X},{\cal E}_{|\{x \} \times X})
\cong {\Bbb C}^{\oplus 2}.
\end{equation}
Therefore $X$ is a connected component of $M_H(v_1)$,
where $v_1=2\rk(E_0)v(E_0)-\omega_X$.

Then ${\cal H}_{\cal E}:{\bf D}(X) \to {\bf D}(X)_{op}$ is an equivalence 
of categories.
As a corollary of this fact, we get that $M_H(v_1)=X$.
By our construction of ${\cal E}$, we see that
\begin{equation}
v({\cal H}_{\cal E}(x))=-(x^{\vee}+2v(E_0)^{\vee}\langle x,v(E_0) \rangle).
\end{equation}
If $E_0={\cal O}_X$ and $v(E)=r+c_1+(s/2) \omega_X$, then
$v({\cal H}_{\cal E}(E))=s+c_1+(r/2) \omega_X$.

From now on,
we assume that $X$ is unnodal, i.e.
there is no $(-2)$-curve.
Let $\sigma$ and $f$ be elliptic curves on $X$ such that
$(\sigma,f)=1$.
Then 
\begin{equation}
 H^2(X,{\Bbb Z})_f=({\Bbb Z}\sigma \oplus {\Bbb Z}f) \perp
 E_8(-1)
\end{equation}
where $H^2(X,{\Bbb Z})_f$ is the torsion free quotient of 
$H^2(X,{\Bbb Z})$.

\begin{prop}\label{prop:00/08/17}
We set $G_1:={\cal E}_{|\{x \} \times X}$ and
$G_2:={\cal E}_{|X \times \{x \} }$.
Assume that $\deg_{G_1}(v)=0$ and 
$l(v):=-\langle v,\omega_X \rangle/\rk v_1>0$, 
$a(v):=\langle v,v_1 \rangle/\rk v_1>0$. 
Let $\varepsilon$ be an element of $K(X) \otimes {\Bbb Q}$
such that $v(\varepsilon) \in
v_1^{\perp} \cap \omega_X^{\perp}$,
$|\langle v(\varepsilon)^2 \rangle| \ll 1$ and
$(H,c_1(\varepsilon))=0$. 
Then ${\cal H}_{\cal E}$ induces an isomorphism
\begin{equation}
{\cal M}_H^{G_1+\varepsilon}(v)^{ss} \to 
{\cal M}_{\widehat{H}}^{G_2+\widehat{\varepsilon}}
(-{\cal H}_{\cal E}(v))^{ss}.
\end{equation}
\end{prop}

\begin{proof}
Since $H$ is general with respect to $v(E_0)$,
we see that ${\cal E}_{|\{x \} \times X}$ is $G_1$-twisted
stable. Then we see that 
Lemma \ref{lem:G^2-2} holds.
We next show that Lemma \ref{lem:G^0-2} holds.
We may assume that $E$ is $\mu$-stable. 
If $E^{\vee \vee} \ne E_0, E_0(K_X)$, 
then $\Hom(E,E_0)=\Hom(E,E_0(K_X))=0$.
If $E^{\vee \vee}=E_0, E_0(K_X)$, then
$\Hom(E,{\cal E}_{|\{x \} \times X})=0$ for
$x \in X \setminus \Supp(E^{\vee}/E)$.
Thus Lemma \ref{lem:G^0-2} holds.
Then the same proof of Theorem \ref{thm:00/08/17}
works and we get our claim.
\end{proof} 

\begin{cor}
$M_H(r-(1/2)\omega_X,{\cal O}_X) \cong \Hilb_X^{(r+1)/2}$
for a general $H$ with respect to
$r-(1/2)\omega_X$.
\end{cor}

\begin{prop}\label{prop:hodge-enriques}
Assume that $r,s>0$.
Then
$e(M_H^{\alpha}(r+c_1-(s/2) \omega_X))=e(M_H^{\alpha}(s-c_1-(r/2) \omega_X))$
for a general $(H,\alpha)$,
if $(c_1^2)<0$, i.e, $\langle v^2 \rangle<rs$,
where $v=r+c_1-(s/2) \omega_X$.
In particular, if $r> \langle v^2 \rangle$, then
we get our claim.
\end{prop}
 
\begin{proof}
If $(c_1^2)<0$, then
Hodge index theorem implies that
there is a divisor $H$ such that
$(H,c_1)=0$ and $({H}^2)>0$. By Riemann-Roch theorem,
we may assume that $H$ is effective. 
Since $X$ is unnodal, $H$ is ample.   
If $E_0={\cal O}_X$, then
$v({\cal E}_{|\{x \} \times X})=2$.
Hence $v$ satisfies assumptions of Proposition \ref{prop:00/08/17}.
Then we get an isomorphism 
\begin{equation}
 M_{H}^{{\cal O}_X+\varepsilon}(r+c_1-(s/2) \omega_X) \to
 M_{H}^{{\cal O}_X+\varepsilon}(s-c_1-(r/2) \omega_X),
\end{equation}
where $(H,{\cal O}_X+\varepsilon)$ is general with respect to $v$.
By Proposition \ref{prop:indep}, we get our claim.
\end{proof}

\begin{thm}\label{thm:enriques}
Let $v=r+c_1-(s/2) \omega_X \in H^*(X,{\Bbb Q})$ 
be a primitive Mukai vector
such that $r$ is odd.
Then
\begin{equation}
 e(M_H(v,L))=e(\Hilb_X^{(\langle v^2 \rangle+1)/2})
\end{equation}
for a general $H$, where $L \in \Pic(X)$ satisfies 
$c_1(L)=c_1$.
In particular,
\begin{enumerate}
\item
$M_H(v) \ne \emptyset$ for a general $H$
if and only if $\langle v^2 \rangle \geq -1$.
\item
$M_H(v,L)$ is irreducible for a general $H$.
\end{enumerate}  
\end{thm}   

\begin{proof}
We first assume that $c_1 \in E_8(-1)$.
We set $l=\gcd(r,c_1)$.
Replacing $v$ by $v \exp(\xi_1)$, $\xi_1 \in E_8(-1)$,
we may assume that $c_1/l$ is primitive and $s> \langle v^2 \rangle$.
Since $v$ is primitive,
$\gcd(l,s)=1$.
By Proposition \ref{prop:hodge-enriques}, we get
\begin{equation}
e(M_H(r+c_1-(s/2) \omega_X))= 
e(M_H(s-c_1-(r/2) \omega_X)).
\end{equation}
Replacing $v=r+c_1-(s/2) \omega_X$ by
$v'=s-c_1-(r/2) \omega_X$, we may assume that $r >\langle v^2 \rangle$.
By the same argument as above,
 we may assume that $l=1$ and $c_1$ is primitive.
We set $D=\sigma-\frac{(\eta^2)}{2}f+\eta$,
where $\eta \in E_8(-1)$ satisfies that $2(\eta,c_1)=s-1$.
Then $(D^2)=0$ and
$-s+2(c_1,D)=-1$.
Since $v \exp(D)=r+(c_1+rD)-1/2 \omega_X$,
$e(M_H(v))=e(M_H(r+(c_1+rD)-1/2 \omega_X))$.
Since $r >\langle v^2 \rangle$,
Proposition \ref{prop:hodge-enriques} implies that 
our claim holds for this case.

We shall next treat the general case.
We use induction on $r$.
We set $c_1:=d_1 \sigma+d_2 f+ \xi$.
Replacing $v$ by $v \exp(k \sigma)$,
we may assume that
$0 \leq |d_1|<r/2$.
We first assume that $d_1 \ne 0$.
We note that $(c_1,f)=d_1$.
Replacing $v$ by $v \exp(\eta)$, $\eta \in E_8(-1)$,
we may assume that $s>\langle v^2 \rangle$.
Then by Proposition \ref{prop:hodge-enriques}, 
$e(M_H(v))=e(M_H(s-c_1-(r/2) \omega_X))$ for a general $H$.
We take an integer $k$ such that
$0<r+2d_1k<2|d_1|<r$.
Then $v \exp(kf)=s+(-c_1+skf)-r'/2 \omega_X$,
where $r'=r+2d_1k$.
Since $s>\langle v^2 \rangle$,
Proposition \ref{prop:hodge-enriques},
implies that $e(M_{H}(s+(-c_1+skf)-r'/2 \omega_X))=
e(M_{H}(r'+(c_1-skf)-(s/2) \omega_X))$ for a general $H$.
By induction hypothesis, we get our claim.

If $d_1=0$, then we may assume that 
$0 \leq |d_2|<r/2$.
If $|d_2|>0$, then we can apply the same argument and get our claim.
If $d_1=d_2=0$, then $c_1 \in E_8(-1)$, so we get our claim.
\end{proof}

\begin{cor}\label{cor:excep}
If $\langle v^2 \rangle=-1$, then there is a stable vector bundle
$E_0$ of $v(E_0)=v$ with respect to $H$.
\end{cor}

\begin{rem}
By the proof,
we also get the following:
Let $v$ be a primitive Mukai vector such that
$\rk v$ is odd.
Then 
$e({\cal M}_H(mv)^{ss})=e({\cal M}_H(m(1-(n/2) \omega_X))^{ss})$.
where $n=\langle v^2 \rangle$.
\end{rem}

\subsection{Relation to Montonen-Olive duality}

We fix an Enriques surface $X$ and fix the following data:
\begin{itemize}
\item
$(H^2(X,{\Bbb Z})_f,Q)$ : a lattice with 
$Q(x,y)=-\int_X x \wedge y, x,y \in H^2(X,{\Bbb Z})_f$.
\item
A orthogonal decomposition of
$H^2(X,{\Bbb R})$ as a sum of definite signature:
\begin{equation}
P:H^2(X,{\Bbb Z}) \otimes {\Bbb R} \cong {\Bbb R}^{9,0} \oplus
 {\Bbb R}^{0,1}.
\end{equation}
\item
$P_L(x)=x_L$, $P_R(x)=x_R$ : the projections onto the two factors.
\end{itemize}

For an odd integer $r>0$, 
we define $U(r)$-partition function by
\begin{equation}\label{eq:u(r)}
 Z^r(\tau,x):=\sum_{
 {\scriptstyle
 v \in H^*(X,{\Bbb Q})}
 \atop{\scriptstyle
 \rk v=r}}
  ``\chi(M_H(v))"q^{\frac{\langle v^2 \rangle}{2r}} 
 q^{\frac{1}{2r}Q(c_1(v)_L^2)}
 \overline{q}^{\frac{-1}{2r}Q(c_1(v)_R^2)} e^{Q(c_1(v),x)},
\end{equation} 
where 
$(\tau, x) \in {\Bbb H} \times H^2(X,{\Bbb Z}) \otimes {\Bbb C}$,
${\Bbb H}:=\{\tau \in {\Bbb C}| \Im \tau>0 \}$,
$q:=\exp(2 \pi \sqrt{-1} \tau)$, $e:=\exp(2 \pi \sqrt{-1})$ and  
$``\chi(M_H(v))"$ 
is a kind of Euler characteristic of a nice compactification of 
$M_H(v)$
(Mathematically there is no definition).
Then Montonen-Olive duality for $U(r)$ gauge group (cf. \cite{V-W:1})
says that
\begin{itemize}
\item[($\#$)]
$Z^r(\tau,x)$ transforms like a Jacobi form
of holomorphic/anti-holomorphic weight 
\begin{equation}
(-\chi(X)/2+b_{-}(X)/2,b_{+}(X)/2)=(-3/2,1/2).
\end{equation} 
\end{itemize}
$SU(1)$-partition function $Z_0^1(\tau)$ is given by
$Z^1_0(\tau)=1/{\eta(\tau)^{12}}$ (\cite{Go:1}, \cite{V-W:1}).
Hence 
\begin{equation}
\begin{split}
Z^1(\tau,x) &=2Z^1_0(\tau)
\left(\sum_{c \in H^2(X,{\Bbb Z})_f} q^{\frac{1}{2}Q(c_L^2)}
 \overline{q}^{\frac{-1}{2}Q(c_R^2)} e^{Q(c,x)} \right)\\
 &= \frac{2}{\eta(\tau)^{12}}\Theta(\tau,P,x),
\end{split}
\end{equation}
where $\Theta(\tau,P,x)= \sum_{c \in H^2(X,{\Bbb Z})_f} 
q^{\frac{1}{2}Q(c_L^2)}
 \overline{q}^{\frac{-1}{2}Q(c_R^2)} e^{Q(c,x)} $
 and the factor $2$ comes from the torsion submodule of
$H^2(X,{\Bbb Z})$.
We set 
\begin{equation}
\Gamma:=\left \langle 
\begin{pmatrix}
1 & 2\\
0 & 1 
\end{pmatrix},
\begin{pmatrix}
0 & 1\\
-1 & 0 
\end{pmatrix} 
\right \rangle
\subset SL_2({\Bbb Z}).
\end{equation}
For an odd integer $r>0$, we set
\begin{equation}
N(r):=\left\{
\begin{pmatrix}
a & b\\
c & d 
\end{pmatrix} 
 \left| 
\begin{split}
a,b,c,d \in {\Bbb Z}\\
ad-bc=r,
\end{split} 
\begin{pmatrix}
a & b\\
c & d 
\end{pmatrix} 
\equiv 
\begin{pmatrix}
1 & 2\\
0 & 1 
\end{pmatrix},
\begin{pmatrix}
0 & 1\\
-1 & 0 
\end{pmatrix} 
\mod 2
\right. \right\}.
\end{equation}
Then we have a coset decomposition:
\begin{equation}
N(r)=\coprod_{
{\scriptstyle a,b,d>0}
\atop{\scriptstyle ad=r, 0 \leq b<d}
}\Gamma
\begin{pmatrix}
a & 2b\\
0 & d 
\end{pmatrix}.
\end{equation} 
Since $\Theta(\tau,P,x)$ transforms like a Jacobi form of 
holomorphic/anti-holomorphic
weight $(9/2,1/2)$,
$Z^1(\tau,x)$ transforms like a Jacobi form of holomorphic/anti-holomorphic
weight $(-3/2,1/2)$:
\begin{equation}
Z^1\left(-\frac{1}{\tau},\frac{x_L}{\tau}+\frac{x_R}{\overline{\tau}}\right)
=(-\sqrt{-1}\tau)^{-3/2}(\sqrt{-1}\overline{\tau})^{1/2}
e^{\frac{Q(x_L^2)}{2\tau}}
e^{\frac{Q(x_R^2)}{2\overline{\tau}}}Z^1(\tau,x_L+x_R).
\end{equation}
We conjecture that
 $Z^r(\tau,x)$ is given by Hecke transformation of 
order $r$ of $Z^1(\tau,x)$ (\cite{MNVW:1}):
\begin{equation}\label{eq:U(r)}
Z^r(\tau,x)=\frac{1}{r^2} 
\sum_{
{\scriptstyle
a,b,d \geq 0}
\atop{\scriptstyle
 ad=r,
b<d}}
d Z^1\left(\frac{a \tau+2b}{d},ax \right).
\end{equation}
Hence $Z^r(\tau,x)$ transforms like a Jacobi form of
holomorphic/anti-holomorphic
weight $(-3/2,1/2)$ and index $r$.
Thus ($\#$) holds.

{\it Evidence:}

We shall give an evidence of this conjecture by
using Theorem \ref{thm:enriques}.
For simplicity, we set $X^{[n]}=\Hilb_X^n$.
Then we see that

\begin{equation}
\begin{split}
& \quad \frac{1}{2}
\sum_{0 \leq b<d}d Z^1\left(\frac{a \tau+2b}{d},ax \right)\\ 
&=
\sum_{0 \leq b<d}\sum_{\xi \in H^2(X,{\Bbb Z})_f}
\sum_n d\chi(X^{[n]})q^{\frac{a}{d}(n-1/2)}
 q^{\frac{a}{2d}Q(\xi_L^2)}
 \overline{q}^{\frac{-a}{2d}Q(\xi_R^2)} e^{aQ(\xi,x)} 
 e^{\frac{2b}{d}((n-1/2)+Q(\xi^2)/2)}\\
 &=\sum_{\xi \in H^2(X,{\Bbb Z})_f}
\sum_{d|2n-1+Q(\xi^2)}
 d^2\chi(X^{[n]})q^{\frac{a}{d}(n-1/2)}
 q^{\frac{a}{2d}Q(\xi_L^2)}
 \overline{q}^{\frac{-a}{2d}Q(\xi_R^2)} e^{aQ(\xi,x)} \\
 &=\sum_{\xi \in H^2(X,{\Bbb Z})_f}
\sum_{k}
 d^2\chi(X^{[(kd+1)/2-Q(\xi^2)/2]})q^{\frac{a}{d}(kd/2-Q(\xi^2)/2)}
 q^{\frac{a}{2d}Q(\xi_L^2)}
 \overline{q}^{\frac{-a}{2d}Q(\xi_R^2)} e^{aQ(\xi,x)} \\
&= \sum_{\xi \in H^2(X,{\Bbb Z})_f}
\sum_{w=(d,\xi,-k/2)}
 d^2\chi(X^{[(\langle w^2 \rangle+1)/2]})
 q^{\frac{a}{2d}\langle w^2 \rangle}
 q^{\frac{a}{2d}Q(c_1(w)_L^2)}
 \overline{q}^{\frac{-a}{2d}Q(c_1(w)_R^2)} e^{aQ(c_1(w),x)} \\
&= \sum_{\rk w=d}
 d^2\chi(X^{[(\langle w^2 \rangle+1)/2]})
 q^{\frac{1}{2r}\langle (aw)^2 \rangle}
 q^{\frac{1}{2r}Q(c_1(aw)_L^2)}
 \overline{q}^{\frac{-1}{2r}Q(c_1(aw)_R^2)} e^{Q(c_1(aw),x)}.
 \end{split}
 \end{equation}
Therefore we get 
\begin{equation}
 ``\chi(M_H(v))"=\sum_{v=aw}\frac{2}{a^2}
 \chi(X^{[(\langle w^2 \rangle+1)/2]}).
 \end{equation}
If $v$ is primitive, by Theorem \ref{thm:enriques}, we get

 \begin{equation}
 ``\chi(M_H(v))"=\chi(X^{[(\langle v^2 \rangle+1)/2]})=\chi(M_H(v)).
 \end{equation}
 This implies that $\chi(M_H(v))$ is related to modular forms  
 and in particular Hecke transforms.

\section{Twisted stability for dimension 1 sheaves}\label{sect:1-dim}

We shall generalize twisted stability to
purely 1-dimensional sheaves.
Let $X$ be a smooth projective surface as in section \ref{sect:pre}.

\begin{defn}
Let $E$ be a purely 1-dimensional sheaf on $X$.
\begin{enumerate}
\item
For a ${\Bbb Q}$-divisor $\alpha$ on $X$ (or its numerical equivalence class),
we define $\alpha$-twisted slope of $E$ by
\begin{equation}
\mu_{\alpha}(E):=\frac{\chi(E(-\alpha))}{(c_1(E),H)}
=\frac{\chi(E)-(c_1(E),\alpha)}{(c_1(E),H)}.
\end{equation}
\item
$E$ is $\alpha$-twisted semi-stable (resp. stable) if
\begin{equation}
\mu_{\alpha}(F) \underset{(<)}{\leq} \mu_{\alpha}(E)
\end{equation}
for all $0 \subsetneq F \subsetneq E$.
\item
For a vector bundle $G$ on $X$ (or its class in $K(X)$), 
we also define $\mu_G(E)$ by
\begin{equation}
\mu_G(E):=\frac{\chi_G(E)}{\deg_G(E)}=\frac{\chi(E)-(c_1(E),c_1(G))/\rk G}{
(c_1(E),H)}=\mu_{c_1(G)/\rk G}(E).
\end{equation}
We also define $G$-twisted stability as $c_1(G)/\rk G$-twisted stability.
\end{enumerate}
\end{defn}

\begin{defn}
Let ${\cal M}_H^{\alpha}(\gamma)^{ss}$ be the moduli stack of
$\alpha$-twisted semi-stable sheaves $E$ of $\gamma(E)=\gamma$
and ${\cal M}_H^{\alpha}(\gamma)^{s}$ the open substack consisting 
of $\alpha$-twisted stable sheaves.
\end{defn}

We next generalize wall and chamber structure.
Let ${\cal D}_{\xi}$ be a set of effective divisors
$D$ such that there is an effective divisor $D'$ and
$D+D' \equiv \xi$.
It is known that ${\cal D}_{\xi}$ is a bounded set.
For $E \in K(X)$ of $\gamma(E)=(0,\xi,\chi)$,
we consider the set ${\cal W}$ of $\alpha \in \NS(X) \otimes {\Bbb Q}$
which satisfies
\begin{equation}
\frac{\chi-(\xi,\alpha)}{(\xi,H)}=\frac{n-(D,\alpha)}{(D,H)},
\end{equation}
for some $n \in {\Bbb Z}$, $D \in {\cal D}_{\xi}$
and $(D,n) \not \in {\Bbb Q}(\xi,\chi)$. 
This set is a countable union of hyperplanes of $\NS(X) \otimes {\Bbb Q}$.
We claim that the numbers of hyperplanes is locally finite.
Proof of the claim: 
Assume that $\alpha$ belongs to a bounded subset of 
$\Pic(X) \otimes {\Bbb Q}$.
Since the choice of $D$ is finite,
the set of $n$ is bounded, which implies that
the choice of $n$ is finite.

\begin{defn}
We call a defining hyperplane of ${\cal W}$ by wall and
a connected component of $\NS(X)\otimes{\Bbb R} \setminus {\cal W}$
by chamber.
\end{defn}

For $\alpha,\alpha'$ which belong to a chamber ${\cal C}$,
${\cal M}_H^{\alpha}(\gamma)^{ss}={\cal M}_H^{\alpha'}(\gamma)^{ss}$.
Hence we may denote this stack by
${\cal M}_H^{\cal C}(\gamma)^{ss}$.
 
\begin{defn}
Let $W$ be a wall and ${\cal C}$ a chamber such that $\overline{{\cal C}}$ intersects $W$. 
Let $\alpha$ be a ${\Bbb Q}$-divisor 
belonging to $\overline{{\cal C}} \cap W$ and $\alpha_1$ 
a ${\Bbb Q}$-divisor which belongs to ${\cal C}$.
$V^{\alpha,{\cal C}}({\gamma})$ be the set of $\alpha$-twisted 
semi-stable sheaves with respect to $H$ 
such that $E$ is not $\alpha_1$-twisted 
semi-stable with respect to $H$ and $\gamma(E)=\gamma$.
\end{defn}

In the same way as in \cite[Prop. 2.4]{Y:2}, we get the following.
\begin{prop}\label{prop:2-2}
Let $C$ be a 2-dimensional vector space such that $C \cap {\cal C} \ne \emptyset$ and $\alpha \in C$.

(1)There is an element $\alpha_1 \in C$ such that
$V^{\alpha,{\cal C}}({\gamma})$ is the set of torsion free sheaves $E$ whose
 Harder-Narasimhan filtration 
\begin{equation}\label{eq:HNF}
0 \subset F_1 \subset F_2 \subset \cdots \subset F_s=E
\end{equation}
with respect to $\alpha_1$
satisfies $\mu_{\alpha}(F_i)=\mu_{\alpha}(F_{i+1})$ and 
is also the Harder-Narasimhan filtration of $E$ with respect to 
$\alpha_t:=\alpha_1+t(\alpha-\alpha_1), 0 \leq t<1$.

(2) ${\cal M}_H^{\alpha}({\gamma})^{ss}={\cal M}_H^{{\cal C}}({\gamma})^{ss}
 \amalg V_H^{\alpha,{\cal C}}({\gamma})$.

\end{prop}

\begin{lem}\label{lem:vanish}
Assume that $-K_X$ is base point free.
Let $E$ and $F$ be $\alpha$-twisted semi-stable sheaves of
$\mu_{\alpha}(E)> \mu_{\alpha}(F)$, or 
of $\mu_{\alpha}(E)= \mu_{\alpha}(F)$ and $(c_1(F),-K_X)>0$.
Then $\Hom(E,F(K_X))=0$.
\end{lem}
\begin{proof}
Considering Jordan-H\"{o}lder filtration,
we may assume that $E$ and $F$ are $\alpha$-twisted stable.
Since $-K_X$ is base point free, there is an inclusion
$F(K_X) \hookrightarrow F$.
Hence we get an inclusion $\Hom(E,F(K_X)) \hookrightarrow
\Hom(E,F)$.
If $\mu_{\alpha}(E)> \mu_{\alpha}(F)$, then 
$\Hom(E,F)=0$, and hence $\Hom(E,F(K_X))=0$. 
Assume that $\mu_{\alpha}(E)= \mu_{\alpha}(F)$.
If $\Hom(E,F) \ne 0$, then $E \cong F$ and
$\Hom(E,F) \cong {\Bbb C}$.
If there is a map $E \to F(K_X)$, then 
$E \to F(K_X) \subset F$ is not an isomorphism, because of
$(c_1(F),-K_X)>0$.
Hence $\Hom(E,F(K_X))=0$.
\end{proof}
\begin{cor}
If $-K_X$ is base point free and $(c_1(x),-K_X)>0$, then
${\cal M}_H^{\alpha}(\gamma(x))^{ss}$ is smooth.
\end{cor}

\begin{prop}\label{prop:wall-cross}
Let $C$ be a 2-dimensional vector space such that 
$C \cap {\cal C} \ne \emptyset$ and $\alpha \in C$.
Assume that $-K_X$ is base point free. Then
\begin{equation}
e({\cal M}_H^{\alpha}(\gamma)^{ss})=
e({\cal M}_H^{\cal C}(\gamma)^{ss})+\sum_{(\gamma_1,\gamma_2,\dots,\gamma_s)}
t^{-\sum_{i<j}(c_1(\gamma_i),c_1(\gamma_j))}
\prod_{i=1}^s e({\cal M}_H^{\cal C}(\gamma_i)^{ss}),
\end{equation}
where $(\gamma_1,\gamma_2,\dots,\gamma_s)$ satisfy
$\mu_{\alpha}(\gamma_i)=\mu_{\alpha}(\gamma_j)$ and
$\mu_{\alpha'}(\gamma_i)>\mu_{\alpha'}(\gamma_j)$ for $i<j$ and
$\alpha' \in C \cap {\cal C}$.
\end{prop}

\begin{proof}
For the Harder-Narasimhan filtration \eqref{eq:HNF},
Lemma \ref{lem:vanish} implies that
\begin{equation}
\Ext^2(F_i/F_{i-1},F_j/F_{j-1})=
\Hom(F_j/F_{j-1},F_i/F_{i-1}(K_X))^{\vee}=0
\end{equation}
 for $j \leq i$.
Moreover we get 
\begin{equation}
\chi(F_i/F_{i-1},F_j/F_{j-1})=
\chi(F_j/F_{j-1},F_i/F_{i-1})=-(c_1(F_i/F_{i-1}),c_1(F_j/F_{j-1})).
\end{equation}
Therefore in the same way as in the proof of \cite[Thm. 3.2]{Y:2}
(see the description of the stack of filtration in
\cite[sect. 5.2]{Y:8}),
we get this proposition.
\end{proof}

By induction on $(H,c_1(\gamma))$ and Proposition \ref{prop:wall-cross}, 
we can show the following proposition.
\begin{prop}\label{prop:e-polynomial}
Assume that $-K_X$ is base point free or $K_X$ is numerically trivial.
Then virtual Hodge polynomial 
$e({\cal M}_H^{\alpha}(\gamma)^{ss})$ does not depend on the choice of
general $\alpha$.
\end{prop}

As in \cite{E-G:1}, we can construct moduli of $\alpha$-twisted semi-stable
sheaves (see Appendix Theorem \ref{thm:alpha-ss}).
Hence we get the following theorem.
\begin{thm}\label{thm:M-W2}
\begin{enumerate}
\item
There is a coarse moduli scheme $\overline{M}_H^{\alpha}(\gamma)$
of $S$-equivalence classes of $\alpha$-twisted semi-stable sheaves.
\item $\overline{M}_H^{\alpha}(\gamma)$ is projective.
\item For different $\alpha,\alpha'$, the relation between 
$\overline{M}_H^{\alpha}(\gamma)$ and $\overline{M}_H^{\alpha'}(\gamma)$
is described as Mumford-Thaddeus type flips:
\begin{equation}
\begin{matrix}
          \overline{M}_H^{\alpha_1}(\gamma) 
           &&&& \overline{M}_H^{\alpha_2}(\gamma)&&&&&
    \overline{M}_H^{\alpha_n}(\gamma)\cr
          &\searrow &&\swarrow&&\searrow&&\cdots & \swarrow &\cr
          &&\overline{M}_H^{\alpha_{1,2}}(\gamma)&&&&
      \overline{M}_H^{\alpha_{2,3}}(\gamma)&&\cr    
\end{matrix}
\end{equation}
where $\alpha=\alpha_1$, $\alpha'=\alpha_n$.
\end{enumerate}
\end{thm}

\section{Fourier-Mukai transform on elliptic surfaces}\label{sect:elliptic}

Let $\pi:X \to C$ be an elliptic surface with a $0$-section $\sigma$  
such that every fiber is irreducible.
Let $f$ be a fiber of $\pi$.
Then a compactification of 
the relative Jacobian of $\pi:X \to C$ is isomorphic to
$\pi:X \to C$.
Let ${\cal P}$ be a universal family on $X \times_C X$.
Then $\widehat{\cal F}_{\cal P}[2]$ is the inverse of ${\cal F}_{\cal P}$.

Let $\tau$ be a section of $\pi$.
We set
\begin{equation}
\langle \sigma,f \rangle^{\perp}:=
\{D \in \NS(X)|(D,\sigma)=(D,f)=0 \}.
\end{equation}
We shall normalize ${\cal P}$ so that
${\cal F}_{\cal P}^1({\cal O}_X)={\cal O}_{\sigma}$.
In \cite[3.2]{Y:7}, we showed the following:
 
For a coherent sheaf $E$ of 
$(\rk(E),c_1(E),-\ch_2(E))=(r,lf+D,n)$, 
$D \in \langle \sigma,f \rangle^{\perp}$
\begin{equation}\label{eq:chern}
 (\rk({\cal F}_{\cal P}(E)),c_1({\cal F}_{\cal P}(E)),\chi({\cal F}_{\cal P}(E)))
 =-(0,r\sigma+nf-D,r+l).
\end{equation}

Let $M(r,d) \to C$ be the relative moduli space of 
stable sheaves of rank $r$ and degree $d$.
\begin{lem}
Assume that $(r,d)=1$. Then
$M(r,d) \cong X$.
\end{lem}
\begin{proof}
We shall prove our claim by induction on $r$.
By tensoring ${\cal O}_X(k\sigma)$, we may assume that
$0< d \leq r$.
If $r=d$, then $r=d=1$. Hence our claim obviously holds.
If $d<r$, then  by $\widehat{{\cal F}}_{\cal P}$, 
we get an isomorphism $M(r,d) \to M(d,-r)$.
By induction hypothesis, $M(d,-r) \cong X$.
Thus our claim holds.
\end{proof}
By this proof, a universal family ${\cal E}$ on
$M(r,d) \times_C X$ is obtained by compositions of $T_{\sigma}$ and
${{\cal F}}_{\cal P}$ from  
${\cal P}$ on $X \times_C X$.
We consider Fourier-Mukai transform 
${\cal F}_{\cal E}:{\bf D}(X) \to {\bf D}(X)$ defined by ${\cal E}$:
\begin{equation}
 {\cal F}_{\cal E}(x)={\bf R}p_{2*}({\cal E} \otimes p_1^*(x)), 
 x \in {\bf D}(X)
\end{equation}
where $p_i:X \times X \to X$, $i=1,2$ are two projections.
\begin{lem}
 \begin{equation}
 {\cal F}_{\cal E}(\Delta)=\Delta
\end{equation}
for $\Delta \in \gamma^{-1}(\langle \sigma,f \rangle^{\perp})$. 
\end{lem}

\begin{proof}
Since $T_{\sigma}(D)=(0,D,(\sigma,D))=(0,D,0)$.
By \eqref{eq:chern}, we get ${\cal F}_{\cal P}(\Delta)=\Delta$.
Since ${\cal E}$ on
$M(r,d) \times_C X$ is obtained by compositions of $T_{\sigma}$ and
$\widehat{{\cal F}}_{\cal P}$ from ${\cal P}$,
Lemma \ref{lem:composition} implies our claim.
\end{proof}

We set $E_0:=\widehat{\cal F}_{\cal E}^1({\cal O}_{\sigma})$
and $F_0:={\cal F}_{\cal E}^0({\cal O}_{\sigma})$.
Then $E_0$ is a vector bundle of $(\rk (E_0),c_1(E_0))=(r,-d \sigma+kf)$,
$k \in {\Bbb Z}$ and
$F_0$ a vector bundle of $(\rk (F_0),c_1(F_0))=(r,d' \sigma+k'f)$,
$k' \in {\Bbb Z}$.
In $K(X)$, we see that
\begin{equation}
 \begin{cases}
  {\cal F}_{\cal E}({\Bbb C})=(F_0)_{|f},\\
  {\cal F}_{\cal E}((E_0)_{|f})=-{\Bbb C},
 \end{cases}
\end{equation}
where ${\Bbb C}$ is the structure sheaf of $f \cap \sigma$.
We see that
\begin{equation}
 K(X)/\gamma^{-1}(\langle \sigma,f \rangle^{\perp}) \subset
 \frac{1}{r}({\Bbb Z}E_0+{\Bbb Z}{E_0}_{|f}+{\Bbb Z}{\Bbb C}).
\end{equation}
For $x=a E_0+b{E_0}_{|f}+c{\Bbb C}$,
${\cal F}_{\cal E}(x)=-a{\cal O}_{\sigma}-b {\Bbb C}+c F_{0|f}$.
Hence
\begin{equation}\label{eq:FM(gamma)}
 -\gamma({\cal F}_{\cal E}(x))=
 (0,a \sigma-crf,b-c \chi({F_0}_{|f})+a \chi({\cal O}_{\sigma})).
\end{equation}
Since $(r,d)=1$, if $x$ belongs to $K(X)$, then $a \in {\Bbb Z}$. 

\begin{lem}\label{lem:f}
If $c_1=\lambda f$, $\lambda \in {\Bbb Z}$, 
then ${\cal M}_{\sigma+nf}^{\alpha}(0,c_1,\chi)^{ss}=
{\cal M}_{\sigma+nf}(0,c_1,\chi)^{ss}$ for any ${\Bbb Q}$-divisor $\alpha$.
In particular, ${\cal M}_{\sigma+nf}^{F_0}(0,rf,d')^{ss}=
{\cal M}_{\sigma+nf}(0,rf,d')^{ss}$.
\end{lem}

\subsection{Preservation of stability}

The following is an easy consequence of the proof of base change theorem
(cf. \cite[Lem. 3.6]{Y:7}).
\begin{lem}\label{lem:base}
Let $L$ be a coherent sheaf of pure dimension 1 on 
$X$ with $c_1(L)=r\sigma+nf+D$, $D \in \langle \sigma,f \rangle^{\perp}$.
Then $L$ satisfies $\WIT_1$ and $\widehat{{\cal F}}_{\cal E}^1(L)$ is torsion free, 
if the following two
conditions are satisfied:
\begin{enumerate}
\item[(1)]
 $\Hom({\cal E}_x,L)=0$, $x \in X$ except finite subset $S$ of $X$.
\item[(2)]
 $\Ext^2({\cal E}_x,L)\cong \Hom(L,{\cal E}_x)^{\vee}=0$ for all $x \in X$. 
\end{enumerate}
\end{lem}

\begin{lem}\label{lem:H0}
For a purely 1-dimensional sheaf $L$ of 
$c_1(L)=r\sigma+nf+D$, $D \in \langle \sigma,f \rangle^{\perp}$,
$\Hom_{p_1}({\cal E},p_2^* L)=0$.
\end{lem}

\begin{proof}
It is sufficient to prove that
$\Hom({\cal E}_x,L)=0$ for some point $x \in X$.
We choose a point $x \in X$ which is not contained in
$\Supp(L)$.
Since $L$ is of pure dimension 1, we get 
$\Hom({\cal E}_x,L)=0$.
\end{proof}

\begin{prop}\label{prop:base-conditions}
Let $L$ be a coherent sheaf of pure dimension 1 on 
$X$ with $c_1(L)=x\sigma+nf+D$, $D \in \langle \sigma,f \rangle^{\perp}$ and
$\chi_{F_0}(L)>0$.
If $L$ is $F_0$-twisted 
semi-stable with respect to $\sigma+kf$, $k \gg 0$, then
$L$ satisfies two conditions in Lemma \ref{lem:base}.
\end{prop}

\begin{proof}
By taking account of a Jordan-H\"{o}lder filtration, we may assume that $L$
is stable.
We note that 
\begin{equation}
 \mu_{F_0}(L):=\frac{\chi_{F_0}(L)}{(c_1(L \otimes F_0^{\vee}),\sigma+kf)}=
 \frac{\chi_{F_0}(L)}{r((k+(\sigma^2))x+n)}.
\end{equation} 
Since $\mu_{F_0}({\cal E}_x)=0$, $x \in X$,
by Lemma \ref{lem:f}, (2) holds.
So we shall prove (1).
Let $D=D_{vir}+D_{hol}$ be the decomposition of the scheme-theoretic
support of $L$, where $D_{vir}$ consists of all fiber components and
$D_{hol}$ consists of the other components.
Then we have an exact sequence
\begin{equation}
 0 \to F \to L \to (L_{|D_{hol}})/T \to 0,
\end{equation}
where $T$ is the torsion submodule of $L_{|D_{hol}}$.   
Then $F$ is a pure dimension 1 subsheaf of $L$ with $c_1(F)=lf$.
By the stability of $L$, we get 
\begin{equation}
 \mu_{F_0}(F)=\frac{\chi_{F_0}(F)}{rl}\leq 
\frac{\chi_{F_0}(L)}{r((k+(\sigma^2))x+n)}.
\end{equation}
Since $k$ is sufficiently large (the condition
$k > \max\{((c_1(L),\sigma+k_0 f)\chi_{F_0}(L)-n)/x-(\sigma^2), k_0\}$ 
is sufficient, where $\sigma+k_0 f$ is ample),
we get $\chi_{F_0}(F) \leq 0$.
Since $\Hom({\cal E}_x,L_{|D_{hol}}/T)=0$ for all $x \in X$,
we shall prove that
$\Hom({\cal E}_x,F)=0$ except finite numbers of points.
\newline
Proof of the claim:
Let 
\begin{equation}
 0 \subset F_1 \subset F_2 \subset \cdots \subset F_s=F
\end{equation}
be the Harder-Narasimhan filtration of $F$ with respect to $\sigma+kf$.
Then 
\begin{equation}
 \mu_{F_0}(F_1)>\mu_{F_0}(F_2/F_1)>\dots>\mu_{F_0}(F_s/F_{s-1}).
\end{equation}
Since $F_1$ is a subsheaf of $L$, 
we also have the inequality $\chi_{F_0}(F_1) \leq 0$. 
If $\Hom({\cal E}_x,F_1) \ne 0$, then 
$\mu_{F_0}(F_1)=0$ and $F_1$ is $S$-equivalent to
${\cal E}_x \oplus E$ for some $E$.
Hence the choice of $x$ is finite.
Clearly $\Hom({\cal E}_x,F_i/F_{i-1})=0$ for $i \geq 2$.
Hence the claim holds. 
\end{proof}

\begin{lem}\label{lem:wit F}
Let $E$ be a torsion free sheaf of $\rk(E)=xr>0$ and
$c_1(E \otimes E_0^{\vee})=lf+D$, 
$D \in \langle \sigma,f \rangle^{\perp}$ on $X$.
Assume that $E$ is $\mu$-semi-stable with respect to $\sigma+k f$, $k \gg 0$.
Then $E$ satisfies $\WIT_1$ and ${\cal F}_{\cal E}^1(E)$ 
is of pure dimension 1.
\end{lem}

\begin{proof}
We shall first prove that ${\cal E} \otimes p_1^*(E)$ is $p_2$-flat.
Let 
\begin{equation}
 0 \to W_1 \to W_0 \to {\cal E} \to 0
\end{equation}
be a locally free resolution of ${\cal E}$ on $X \times X$.
It is sufficient to prove that
\begin{equation}
 \psi_x:(W_1)_{|x \times X} \otimes E \to (W_0)_{|x \times X}\otimes E
\end{equation}
is injective for all $x \in X$.  
We note that $\rk W_1=\rk W_0$ and ${\cal E}_x \otimes E$ is a 
torsion sheaf on $X$.
Since $E$ is torsion free, $\psi_x$ is injective for all $x \in X$.
Thus ${\cal E} \otimes p_1^*(E)$ is a $p_2$-flat sheaf.

Hence we can use base change theorem.
Since $p_2:X \times_C X \to X$ is relative dimension 1,
$R^2p_{2*}({\cal E} \otimes p_1^*(E))=0$.
Since $E_{|\pi^{-1}(y)}$ is semi-stable for general $y \in C$,
$H^0(X,{\cal E}_x \otimes E)=0$ for a general point $x$ of $X$.
Thus $p_{2*}({\cal E} \otimes p_1^*(E))$ is a torsion sheaf.
By the proof of base change theorem,
locally there is a complex of locally free sheaves 
$V_1 \overset{\phi}\to V_0$ which is quasi-isomorphic to
${\mathbf R}p_{2*}({\cal E}\otimes p_1^*(E))$.
Hence $p_{2*}({\cal E} \otimes p_1^*(E))=\ker \phi=0$,
which means that $E$ satisfies $\WIT_1$.
Also we get ${\mathrm {proj-}}\dim(\coker \phi)=1$.
Hence $R^1p_{2*}({\cal E} \otimes p_1^*(E))$ is of pure dimension 1.
\end{proof}

\begin{cor}\label{cor:cal F}
Let $E$ be a torsion free sheaf on $X$ and assume that
$E_{|\pi^{-1}(y)}$ is a semi-stable vector bundle of 
$\deg(E(-(d/r)\sigma)_{|\pi^{-1}(y)})=0$
for a general $y \in C$.
Then $E$ satisfies $\WIT_1$ and 
${\cal F}_{\cal E}^1(E)$ is of pure dimension 1. 
\end{cor}

\begin{proof}
We consider the Harder-Narasimhan filtration of $E$
with respect to $\sigma+kf$, $k\gg 0$.
Applying ${\cal F}_{\cal E}$ to this filtration, we get our corollary
by Lemma \ref{lem:wit F}. 
\end{proof}

\begin{lem}\label{lem:cal F}
We set $\gamma(E)=\gamma(xE_0+y {E_0}_{|f}+z {\Bbb C})+D$,
$D \in \langle \sigma,f \rangle^{\perp}$,
and assume that $E$ is $E_0$-twisted semi-stable with respect to 
$\sigma+k f$, $k \gg 0$.
If $\chi_{F_0}({\cal F}_{\cal E}^1(E))=
x(r\chi({\cal O}_{\sigma})-d'(\sigma^2)-k')+ry>0$, 
then ${\cal F}_{\cal E}^1(E)$ is $F_0$-twisted semi-stable.
\end{lem}

\begin{proof}
Assume that ${\cal F}_{\cal E}^1(E)$ is not semi-stable.
Then, there is a stable subsheaf $F$ of ${\cal F}_{\cal E}^1(E)$ such that 
$ \mu_{F_0}(F)>\mu_{F_0}({\cal F}_{\cal E}^1(E))=
(\chi_{F_0}({\cal F}_{\cal E}^1(E))/((k+(\sigma^2))x-rz)>0$.
We set $G:={\cal F}_{\cal E}^1(E)/F$.
Applying $\widehat{\cal F}_{\cal E}$ to the exact sequence
\begin{equation}
 0 \to F \to  {\cal F}_{\cal E}^1(E) \to G \to 0,
\end{equation}
we get an exact sequence
\begin{equation}
 \begin{CD}
  0 @>>>\widehat{\cal F}_{\cal E}^0(F) @>>> 0 @>>>\widehat{\cal F}_{\cal E}^0(G)\\
  @>>> \widehat{\cal F}_{\cal E}^1(F) @>>> E @>>>\widehat{\cal F}_{\cal E}^1(G)\\ 
  @>>> \widehat{\cal F}_{\cal E}^2(F) @>>> 0 @>>>\widehat{\cal F}_{\cal E}^2(G) @>>> 0 
 \end{CD}
\end{equation}
By Lemma \ref{lem:H0}, $\widehat{\cal F}_{\cal E}^0(G)=0$.
Since $\mu_{F_0}(F)>0$, we also get $\widehat{\cal F}_{\cal E}^2(F)=0$.
Thus $F$ and $G$ satisfies $\WIT_1$ and $\widehat{\cal F}_{\cal E}^1(F)$ 
is a subsheaf of $E$.
We set
$\gamma(\widehat{\cal F}_{\cal E}^1(F))=
\gamma(x'E_0+y'{E_0}_{|f}+z'{\Bbb C})+D'$,
$D' \in \langle \sigma,f \rangle^{\perp}$.
Then $\gamma(F)=
(0,x' \sigma-z'rf-D',y'-z'\chi({F_0}_{|f})+x'\chi({\cal O}_{\sigma}))$.
Since $E$ is semi-stable with respect to $\sigma+kf$, $k \gg 0$,
\begin{itemize}
\item[($\star$)]
\begin{enumerate}
 \item
$y'/x'<y/x$,
or
\item
$y'/x'=y/x$ and 
  $z'/x'\leq z/x$.
\end{enumerate}
\end{itemize}
On the other hand,
\begin{equation}
 \begin{split}
  \mu_{F_0}({\cal F}_{\cal E}^1(E))-\mu_{F_0}(F)=&
  \frac{xl+yr}{r(k+(\sigma^2))x-rz}-
  \frac{x'l+y'r}{r(k+(\sigma^2))x'-rz'}\\
  =&\frac{(yx'-y'x)r(k+(\sigma^2))+
  (x'l+ry')z-(xl+ry)z'}
  {(r(k+(\sigma^2))x-rz)((k+(\sigma^2))x'-z')}\\
  =&\frac{(yx'-y'x)(r(k+(\sigma^2))-rz/x)+
  (xl+ry)(x'z/x-z')}
  {(r(k+(\sigma^2))x-rz)((k+(\sigma^2))x'-z')},
 \end{split}
\end{equation}
where $l=r\chi({\cal O}_{\sigma})-d'(\sigma^2)-k'$.
We note that the choice of $x'$ is finite.
In Lemma \ref{lem:bdd},
we shall show that the choice of $z'$ is also finite.
Then there is an integer $N(x,D,y,z)$ such that for $k>N(x,D,y,z)$,
\begin{enumerate}
\item 
$yx'-y'x>0$ implies 
$(yx'-y'x)(r(k+(\sigma^2))-rz/x)+(xl+ry)(x'z/x-z')>0$
 and 
\item
$yx'-y'x<0$ implies 
$(yx'-y'x)(r(k+(\sigma^2))-rz/x)+(xl+ry)(x'z/x-z')<0$.
\end{enumerate}
Then $(\star)$ implies that 
$\mu_{F_0}({\cal F}_{\cal E}^1(E))-\mu_{F_0}(F) \geq 0$,
which is a contradiction.
Therefore ${\cal F}_{\cal E}^1(E)$ is a semi-stable sheaf.
\end{proof}

\begin{lem}\label{lem:bdd}
Keep the notations as above.
Then the choice of $z'$ is finite and 
the number of such $z'$ is bounded in terms of
$(x,D,y,z)$.
\end{lem}
\begin{proof}
We fix an ample divisor $\sigma+k_0f$.
Since $F$ is a subsheaf of ${\cal F}_{\cal E}^1(E)$,
\begin{equation}
 0 \leq (c_1(F),\sigma+k_0f) \leq (c_1({\cal F}_{\cal E}^1(E)),\sigma+k_0f).
\end{equation}
Since $(c_1(F),\sigma+k_0f)=(k_0+(\sigma^2))x'-rz'$ and 
$x' \leq x$, we get our claim.
\end{proof}

\begin{lem}
Let $L$ be a pure dimension 1 sheaf of 
$\gamma(L)=(0,x\sigma-rzf,y-z\chi({F_0}_{|f})+x\chi({\cal O}_{\sigma}))$.
Assume that $\chi_{F_0}(L)>0$ and 
$L$ is $F_0$-twisted semi-stable with respect to $\sigma+kf$, $k\gg0$.
Then $\widehat{\cal F}_{\cal E}^1(L)$ is an $E_0$-twisted semi-stable sheaf
with respect to $\sigma+kf$, $k\gg0$.
\end{lem}

\begin{proof}
We note that Lemma \ref{lem:base} and Proposition
\ref{prop:base-conditions} imply that 
$L$ satisfies $\WIT_1$ and $\widehat{\cal F}_{\cal E}^1(L)$ is torsion free.
Assume that $\widehat{\cal F}_{\cal E}^1(L)$ 
is not semi-stable with respect to $\sigma+kf$, $k\gg0$.
Then there is a destabilizing subsheaf $F$ of $\widehat{\cal F}_{\cal E}^1(L)$.
We set $G:=\widehat{\cal F}_{\cal E}^1(L)/F$.
It is easy to see that $\widehat{\cal F}_{\cal E}^1(L)_{|\pi^{-1}(y)}$ is semi-stable for
general $y \in C$.
Since $k$ is sufficiently large,
$F_{|\pi^{-1}(y)}$ and $G_{|\pi^{-1}(y)}$ are semi-stable vector bundles
of degree 0 for general $y \in C$.
Then Corollary \ref{cor:cal F} implies that
 $F$ and $G$ satisfy $\WIT_1$ and we get an exact sequence
\begin{equation}
 0 \to {\cal F}_{\cal E}^1(F) \to L \to {\cal F}_{\cal E}^1(G) \to 0.
\end{equation}
In the same way as in Lemma \ref{lem:cal F},
we get a contradiction.
Thus $\widehat{\cal F}_{\cal E}^1(L)$ is semi-stable.
\end{proof}

Therefore, we get the following theorem. 
\begin{thm}\label{thm:FM}
Keep the notation as abobe.
For $x \in K(X)$ of $\rk x>0$,
Assume that the relative twisted degree
$(c_1(x \otimes E_0^{\vee}),f)=0$. 
Then ${\cal F}_{\cal E}$ gives an isomorphism of moduli stack
\begin{equation}
{\cal M}^{E_0}_{\sigma+kf}(\gamma(x))^{ss} \to 
{\cal M}^{F_0}_{\sigma+kf}(-\gamma({\cal F}_{\cal E}(x)))^{ss}
\end{equation}
if $\chi_{F_0}(-{\cal F}_{\cal E}(x))>0$ and $k \gg 0$.
\end{thm}

\begin{rem}\label{rem:epsilon}
Let $\varepsilon \in \langle \sigma,f \rangle^{\perp}$ 
be a ${\Bbb Q}$-divisor of $-(\varepsilon^2) \ll 1$.
Then we can show that
\begin{equation}
{\cal M}^{E_0}_{\sigma+kf+\varepsilon}(\gamma(x))^{ss} \to 
{\cal M}^{F_0(-\varepsilon)}_{\sigma+kf}(-\gamma({\cal F}_{\cal E}(x)))^{ss}
\end{equation}
if $(c_1(x \otimes E_0^{\vee}),f)=0$,
$\chi_{F_0}(-{\cal F}_{\cal E}(x))>0$ and $k \gg 0$. 
\end{rem}

\subsection{Application of Theorem \ref{thm:FM}}
 
\begin{lem}
Let $H=\sigma+kf$ be an ample divisor.
We set $\gamma_n:=(0,\sigma+nf+D,\chi)$, $D \in \langle \sigma,f \rangle$.
Then for a general $t \in {\Bbb Q}$,
${\cal M}_H^{t(\sigma+mf)}(\gamma_n)^{ss}=
{\cal M}_H^{t(\sigma+mf)}(\gamma_n)^s$, if $m \ne k$.
\end{lem}

\begin{prop}\label{prop:e-polynomial2}
We set $\gamma_n:=(0,\sigma+nf+D,\chi)$, 
$D \in \langle \sigma,f \rangle^{\perp}$.
Then $e(M_H^{\alpha}(\gamma_n))$ does not depend on $\alpha$,
if $\alpha$ does not lie on walls.
\end{prop}

\begin{proof}
It is sufficient to prove the claim for
${\cal M}_H^{\alpha}(\gamma_n)^{ss}$.
We prove our claim by induction on $n$.
We note that there is a section $\tau$ of $\pi$ such that
$D=\tau-\sigma-(\tau-\sigma,\sigma)f$.
Hence $\sigma+nf+D=\tau+n'f$, $n'=n-(\tau-\sigma,\sigma)$.
Let $W$ be a wall and $\alpha$ belongs to $W$. 
Let $\alpha_+, \alpha_-$ be ${\Bbb Q}$-divisors
which are very close to $\alpha$ and $\alpha=(\alpha_{+}+\alpha_{-})/2$.
We consider Harder-Narasimhan filtration \eqref{eq:HNF}
in Proposition \ref{prop:2-2}, where $\alpha_1=\alpha_+, \alpha_-$.
Since $c_1(E)=\tau+n'f$, we get
$c_1(F_i/F_{i-1})=l_if, \tau+l_if$.
Assume that $c_1(F_i/F_{i-1})=l_if$ and $c_1(F_j/F_{j-1})=l_jf$
for different $i$ and $j$.
Since $\mu_{\alpha}(F_i/F_{i-1})=\mu_{\alpha}(F_j/F_{j-1})$,
$(F_i/F_{i-1})/l_i=(F_j/F_{j-1})/l_j$ in $K(X) \otimes {\Bbb Q}$.
Then $\mu_{\alpha_{\pm}}(F_i/F_{i-1})=\mu_{\alpha_{\pm}}(F_j/F_{j-1})$,
which is a contradiction.
Therefore $s=2$ and $c_1(F_i/F_{i-1})=l_if$
and $c_1(F_j/F_{j-1})=\tau+l_jf$, $\{i,j \}=\{1,2\}$. 
Since $F_1$ or $F_2/F_1$ is supported on some fibers,
we see that 
\begin{equation}
\Ext^2(F_2/F_1,F_1)=\Hom(F_1,F_2/F_1(K_X))^{\vee}=
\Hom(F_1,F_2/F_1)^{\vee}=0.
\end{equation}
Since $-\chi(F_2/F_1,F_1)=(c_1(F_2/F_1),c_1(F_1))$,
by the same argument as in Proposition \ref{prop:wall-cross},
we get  
\begin{equation}
e({\cal M}_H^{\alpha}(\gamma_n)^{ss})=
e({\cal M}_H^{\alpha_{\pm}}(\gamma_n)^{ss})+
\sum_k e({\cal M}_H^{\alpha_{\pm}}(0,\tau+(n'-kl)f,\chi-kd)^{ss})
e({\cal M}_H^{\alpha_{\pm}}(0,klf,kd)^{ss})t^{kl}
\end{equation}
where $(l,d)$ satisfies that
$\mu_{\alpha}(\gamma_n)=(d-l(\alpha,f))/l(H,f)$.
By induction hypothesis,
\begin{equation}
e({\cal M}_H^{\alpha_+}(0,\tau+(n'-kl)f,\chi-kd)^{ss})=
e({\cal M}_H^{\alpha_-}(0,\tau+(n'-kl)f,\chi-kd)^{ss}).
\end{equation}
Also we know that ${\cal M}_H^{\beta}(0,klf,kd)^{ss})$ 
does not depend on $\beta$.
Therefore we get our claim.
\end{proof}

The following is a generalization of \cite{Go:2} and \cite{Y:number}.

\begin{thm}
If $(r,(c_1,f))=1$, then 
\begin{equation}
e(M_{\sigma+kf}(r,c_1,\chi))
=e(\Pic^0(X) \times \Hilb_X^n),
\end{equation}
where $2n+h^1({\cal O}_X)=\dim M_{\sigma+kf}(r,c_1,\chi)$
and $k \gg 0$.
\end{thm}
\begin{proof}
Let $x$ be an element of $K(X)$
such that $\gamma(x)=(r,c_1,\chi)$.
We set $d:=(c_1,f)$.
We consider Fourier-Mukai transform induced by a universal family
${\cal E}$ on $X \times_C M(r,d)=X \times_C X$.
Then $-c_1({\cal F}_{\cal E}(x))=\sigma+nf+D$,
$n \in {\Bbb Z}$, $D \in \langle \sigma,f \rangle^{\perp}$.
Replacing $x$ by $x \otimes {\cal O}_X(mf)$, $m \gg 0$,
we may assume that $\chi(-{\cal F}_{\cal E}(x))>0$ and
$\chi_{F_0}(-{\cal F}_{\cal E}(x))>0$.
By Theorem \ref{thm:FM},
$M_{\sigma+kf}(\gamma(x))$ is isomorphic to 
$M_{\sigma+kf}^{\alpha}(-\gamma({\cal F}_{\cal E}(x)))$.
By Proposition \ref{prop:e-polynomial2},
$e(M_{\sigma+kf}^{\alpha}(-\gamma({\cal F}_{\cal E}(x))))=
e(M_{\sigma+kf}^{0}(-\gamma({\cal F}_{\cal E}(x))))$.
By using Theorem \ref{thm:FM} again,
we see that
$M_{\sigma+kf}^{0}(-\gamma({\cal F}_{\cal E}(x)))$ 
is isomorphic to Hilbert scheme of points.
\end{proof}
By the similar arguments as in the proof of
Proposition \ref{prop:e-polynomial2},
we also get the following proposition. 
\begin{prop}
$e(M_H^{\alpha}(\gamma_n))$ does not depend on the choice of $H$, 
if $H$ is general.
\end{prop}

\begin{prop}
We set $\gamma_n=(0,\tau+nf,m)$, where $\tau$ is a section of $\pi$.
Then
$M_H^{\alpha}(\gamma_n)$ is smooth for general $\alpha$ and $H$.
\end{prop}

\begin{proof}
Let $F$ be a simple pure dimension 1 sheaf of
$c_1(F)=\tau+nf$.
It is sufficient to show that $H^0(X,K_X) \to \Hom(F,F(K_X))$ is 
an isomorphism.
We set $F_2:=F_{|\sigma}/(torsion)$ and $F_1:=\ker(F \to F_2)$. 
Since $F_1$ is supported on fibers, $F_1(K_X) \cong F_1$.
By the simpleness of $F$,
$\Hom(F,F_1(K_X)) \cong \Hom(F,F_1)=0$.
We note that $\Hom(F_1,F_2(K_X))=0$.
Hence $\Hom(F,F(K_X)) \to \Hom(F,F_2(K_X))$ is injective
and $\Hom(F_2,F_2(K_X)) \to \Hom(F,F_2(K_X))$ is isomorphic.
Since $F_2$ is a line bundle on $\tau$,
\begin{equation}
 H^0(C,\pi_*(K_X))=H^0(X,K_X) \to \Hom(F_2,F_2(K_X))
\end{equation}
is isomorphic.
Therefore $H^0(X,K_X) \to\Hom(F,F(K_X))$ must be isomorphic.
\end{proof}

From now on, we assume that $X$ is a rational elliptic surface.
By using Proposition \ref{prop:e-polynomial} and Remark \ref{rem:epsilon}, 
we also get the following. 
\begin{cor}\label{cor:rational1}
Assume that $X$ is a rational elliptic surface.
Let $\varepsilon, \varepsilon' \in \langle \sigma,f \rangle^{\perp}$
be general ${\Bbb Q}$-divisors of $-(\varepsilon^2),-({\varepsilon'}^2) \ll 1$.
If $\gcd(r,(c_1,f))=l$, then
\begin{equation}
e({\cal M}_{\sigma+kf+\varepsilon}(r,c_1,\chi)^{ss})=
e({\cal M}_{\sigma+kf+\varepsilon'}(l,\xi,\chi')^{ss}), k \gg 0
\end{equation}
for some $\chi' \in {\Bbb Z}$ and
$\xi \in H^2(X,{\Bbb Z})$ of $(\xi,f)=0$.
In particular, if $l=2$, then we obtain 
$e({\cal M}_{\sigma+kf+\varepsilon}(r,c_1,\chi)^{ss})$ from \cite{Y:dual}.
\end{cor}
We also get the following corollary which is a generalization of \cite{Y:dual}.
\begin{cor}
Keep notation as in Corollary \ref{cor:rational1}.
If $(r,D,\chi)$, $D \in \langle \sigma,f \rangle^{\perp}$ 
is primitive, then
\begin{equation}
e({\cal M}_{\sigma+kf+\varepsilon}(r,D,\chi)^{ss})
=e({\cal M}_{\sigma+kf+\varepsilon'}(r,D+tf,\chi)^{ss}), k \gg 0
\end{equation}
for $t \in {\Bbb Z}$.
\end{cor}  

\begin{proof}
If $(r,D,\chi)$, $D \in \langle \sigma,f \rangle^{\perp}$ 
is primitive, then $C:=r \sigma+(r-\chi)f-D$ is primitive.
Since $\NS(X)=H^2(X,{\Bbb Z})$ is a unimodular lattice,
there is a divisor $\beta$ such that $(C,\beta)=1$.
Hence $e({\cal M}_{\sigma+kf}^{\alpha+t \beta}(0,C,m)^{ss})=
e({\cal M}_{\sigma+kf}^{\alpha}(0,C,m-t)^{ss})$, $t \in {\Bbb Z}$.
Therefore $e({\cal M}_{\sigma+kf+\varepsilon}(r,D,\chi)^{ss})
=e({\cal M}_{\sigma+kf+\varepsilon'}(r,D+tf,\chi)^{ss})$.
\end{proof}

\section{Appendix}

Let $X$ be a projective surface defined over a field $k$. 
Let $\alpha$ be a ${\Bbb Q}$-divisor on $X$.
In this appendix, 
we shall prove the following.
\begin{thm}\label{thm:alpha-ss}
There is a projective moduli scheme 
$\overline{M}_H^{\alpha}(\gamma)$ of 
$\alpha$-twisted semi-stable sheaves of pure dimension 1 on $X$
with respect to $H$.
\end{thm}
Since $\alpha$-twisted semi-stability does not change
under the operation $E \mapsto E \otimes {\cal O}_X(nH)$,
we may assume that $\alpha=D/n$ for some positive integer $n$ and
an irreducible divisor $D$ on $X$ of $(D,H)>(c_1,H)$.
We note that 
\begin{equation}
\begin{split}
\chi(E(-\alpha))&=\frac{1}{n}\chi(E(-D))+\left(1-\frac{1}{n}\right)\chi(E)\\
&=\chi(E(-D))+\left(1-\frac{1}{n}\right)\chi(E_{|D}).
\end{split}
\end{equation}
As in \cite{E-G:1},
we can regard $\chi(E(-\alpha))$
as the parabolic Euler characteristic of
the parabolic sheaf $E(-D) \subset E$, where the weight is 
$1-1/n$. 
Then it is sufficient to construct 
moduli space of parabolic semi-stable sheaves
on $X$. 
The moduli space of parabolic stable sheaves 
was constructed by Maruyama and Yokogawa
\cite{M-Y:1} and Inaba \cite{Inaba:1}.
They constructed moduli space as GIT quotient of a suitable space.
Then the problem is to analyse properly semi-stable points of this space.
For a nonsingular variety,
Yokogawa \cite{Yk:1} constructed the moduli space of parabolic semi-stable
torsion free sheaves.  
By technical reason, we can only treat the dimension 1 case.

\subsection{Construction of moduli spaces}

Let $(X,{\cal O}_X(1))$ be a projective scheme $X$ 
over a field $k$ and an ample line bundle 
${\cal O}_X(1)$ an $X$.
Let us start with
the definition of parabolic sheaves.

\begin{defn}
Let $E$ be a pure dimensional $1$ coherent sheaf
on $X$ such that $\dim(D \cap \Supp(E))=0$.
Let 
\begin{equation}
F(E):E(-D)=F_{l+1}(E) \subset F_l(E) \subset \cdots \subset
F_2(E) \subset F_1(E)=E
\end{equation}
be a filtration of coherent sheaves and
$0<\alpha_1 \leq \alpha_2 \leq \dots \leq \alpha_l \leq 1$ 
a sequence of rational numbers.  
Then the triple $E_*:=(E,F(E),\alpha_*)$ is called a parabolic sheaf on $X$,
where $\alpha_*=(\alpha_1,\alpha_2,\dots,\alpha_l)$.
\end{defn}

We shall generalize the notion of parabolic sheaves as follows.
\begin{defn}
Let $E$ be a purely $1$ dimensional coherent sheaf
on $X$.
Let 
\begin{equation}
F(E):F_{l+1}(E) \subset F_l(E) \subset \cdots \subset
F_2(E) \subset F_1(E)=E
\end{equation}
be a filtration of coherent sheaves 
such that $\dim (E/F_{l+1}(E))=0$ and
$0<\alpha_1 \leq \alpha_2 \leq \dots \leq \alpha_l \leq 1$ 
a sequence of rational numbers.  
Then we call the triple $E_*:=(E,F(E),\alpha_*)$ a generalized 
parabolic sheaf on $X$.
\end{defn}
Thus a generalized parabolic sheaf is a sheaf with
a special type of filtration and a sequence of rational numbers.


\begin{defn}
Let $E_*:=(E,F(E),\alpha_*)$ be a generalized parabolic sheaf.
We set
\begin{equation}
\text{par-}\chi(E_*(m)):
=\chi(F_{l+1}(m))+\sum_{i=1}^l \alpha_i \chi(gr_i(E)(m)),
\end{equation}
where $gr_i(E)=F_i(E)/F_{i+1}(E)$.
\end{defn}
We set $\varepsilon_i:=\alpha_{i+1}-\alpha_i$, $1 \leq i \leq l$,
where $\alpha_{l+1}=1$.
Then 
\begin{equation}
\text{par-}\chi(E_*(m))
=\chi(E(m))-\sum_{i=1}^l \varepsilon_i \chi(gr_i(E)(m)).
\end{equation}
For a coherent sheaf $E$ of dimension $1$, we define 
$a_1(E) \in {\Bbb Z}$ by 
\begin{equation}
\chi(E(m))=a_1(E) m+\chi(E).
\end{equation}
For a numerical polynomial of degree 1 
\begin{equation}
h(x)=a_1 x+a_0, a_i \in {\Bbb Z},
\end{equation} 
we set $a_i(h):=a_i$. 

\begin{defn}
A generalized parabolic sheaf $E_*$ of dimension $1$
is sem-stable (resp. stable) with respect to ${\cal O}_X(1)$,
if 
\begin{equation}
\begin{split}
\frac{\text{par-}\chi(E'_*)}{a_1(E')} & \leq 
\frac{\text{par-}\chi(E_*)}{a_1(E)}\\
&(<)
\end{split}
\end{equation}
for all non-trivial generalized parabolic subsheaf $E'_*$ of
$E_*$.
\end{defn}
As in \cite{M-Y:1}, we can consider Jordan-H\"{o}lder filtration
and define $S$-equivalence classes of
generalized parabolic semi-stable sheaves.  

Let $h(x)$ be a numerical polynomial of degree $1$ and $h_i(x)$,
$1 \leq i \leq l$ constant numerical polynomials.
We set $h_*(x):=(h(x),h_1(x),\dots,h_l(x))$.
We shall construct the moduli space $\overline{M}_{X/{\Bbb C}}^{h_*,\alpha_*}$
of $S$-equivalence classes of
generalized parabolic semi-stable sheaves
$(E,F(E),\alpha_*)$ of Hilbert polynomials 
\begin{equation}
 (\chi(E(m)),\chi(gr_1(E)(m)),\dots,
 \chi(gr_l(E)(m)))=h_*(m).
\end{equation}
Let $D$ be a Cartier divisor on $X$.
Then we can construct the 
moduli space $\overline{M}_{D/X/{\Bbb C}}^{h_*,\alpha_*}$
of $S$-equivalence classes of
parabolic semi-stable sheaves as a locally closed subscheme of 
$\overline{M}_{X/{\Bbb C}}^{h_*,\alpha_*}$.

\begin{defn}
Let $\lambda$ be a rational number. 
Then a purely $1$-dimensional sheaf $E$ is of type $\lambda$, if
\begin{equation}
 \frac{\chi(E'')}{a_1(E'')} \geq 
 \frac{\chi(E)}{a_1(E)}-\lambda
\end{equation}
for all quotient sheaf $E''$ of pure dimension $1$.
\end{defn}

\begin{lem}\label{lem:lambda}
Let $E_*=(E,F(E),\alpha_*)$ 
be a generalized parabolic semi-stable sheaf.
Then $E$ is of type $\sum_i \varepsilon_i \chi(gr_i(E))/a_1(E)$.
\end{lem}

\begin{proof}
Let $E \to E''$ be a quotient such that $E''$ is of pure dimension 
$1$. Let $E''_*$ be the induced generalized parabolic structure.
Then 
\begin{equation}
\frac{\chi(E'')-\sum_i \varepsilon_i \chi(gr_i(E''))}{a_1(E'')} \geq 
 \frac{\chi(E)-\sum_i \varepsilon_i \chi(gr_i(E))}{a_1(E)}.
\end{equation}
Hence $\chi(E'')/a_1(E'') \geq \chi(E)/a_1(E)-
\sum_i \varepsilon_i \chi(gr_i(E))/a_1(E)$.
\end{proof}

Let $E_*$ be a generalized parabolic sheaf such that 
$E$ is of type $\lambda$.
Let $E_* \to E_*''$ be a quotient generalized parabolic sheaf such that
$E''$ is of pure dimension $1$ and
\begin{equation}\label{eq:a3}
\frac{\chi(E)-\sum_i \varepsilon_i \chi(gr_i(E))}{a_1(E)} \geq 
\frac{\chi(E'')-\sum_i \varepsilon_i \chi(gr_i(E''))}{a_1(E'')}.
\end{equation}
Then $\chi(E)/a_1(E)+\sum_i \varepsilon_i \chi(gr_i(E))
\geq \chi(E'')/a_1(E'')$.
Since the set of $E_*$ is bounded,
by Grothendieck's boundedness theorem, the set of such quotients
$E''_*$ is bounded.
Hence there is an integer $m({\lambda})$ which depends on $h_*, \alpha_*$
and $\lambda$ such that, for $m \geq m({\lambda})$
and the kernel $E_*'$ of $E_* \to E_*''$ which satisfies \eqref{eq:a3},
\begin{enumerate}
 \item[($\flat 1$)]
  $F_i(E')(m)$, $1 \leq i \leq l+1$ are generated by global sections and
 \item[($\flat 2$)]
  $H^1(X,F_i(E')(m))=0$.
\end{enumerate}
In particular, 
\begin{enumerate}
\item
$E(m)$ is generated by global sections and $H^1(X,E(m))=0$,
\item 
$H^1(X,F_i(E)(m))=0$.
\end{enumerate}

Let $V_m$ be a vector space of dimension $h(m)$.
Let ${\cal Q}:=\Quot_{V_m \otimes {\cal O}_X/X}^{h[m]}$ be the quot-scheme and
$V_m \otimes {\cal O}_{{\cal Q} \times X} \to \widetilde{E}$
the universal quotient sheaf, where $h[m](x)=h(m+x)$.
We set ${\cal Q}_i:=\Quot_{\widetilde{E}/{\cal Q} \times X/{\cal Q}}^{h_i[m]}$
, $1 \leq i \leq l$ and let 
$\widetilde{E} \otimes {\cal O}_{{\cal Q}_i \times X} \to \widetilde{E}_i$
be the universal quotient sheaf.
Then there is a closed subscheme $\Gamma$ of $\prod_{i=1}^l {\cal Q}_i$
which parametrizes sequences of quotients
\begin{equation}\label{eq:a4}
V_m \otimes {\cal O}_X \to E(m) \to E_l(m) \to E_{l-1}(m) \to \dots 
\to E_1(m) \to 0
\end{equation}
such that $V_m \otimes {\cal O}_X \to E(m) \in {\cal Q}$ and 
$E(m) \to E_i(m) \in {\cal Q}_i$, $1 \leq i \leq l$.
We set $F_{i+1}(E):=\ker(E \to E_i)$, $1 \leq i \leq l$ and
$F_1(E):=E$.
Then we have a filtration
\begin{equation}\label{eq:filter}
F(E): F_{l+1}(E) \subset F_{l}(E) \subset \dots \subset
F_1(E)=E.
\end{equation}
Thus $F(E)$ and $\alpha_*$ give a structure of generalized 
parabolic sheaf on $E$.
%

Let $\Gamma^{ss}$ be the open subscheme of $\Gamma$ consisting of
quotients \eqref{eq:a4}
such that 
\begin{enumerate}
\item
$V_m \to H^0(X,E(m))$ is an isomorphism and 
\item
$(E,F(E),\alpha_*)$ is 
a generalized parabolic semi-stable sheaf.
\end{enumerate}

Let $G(n):=Gr(V_m \otimes W,h[m](n))$
be the Grassmannian parametrizing $h[m](n)$-dimensional quotient
spaces of $V_m \otimes W$.
We set $G_i:=Gr(V_m,H_i(m))$.
Assume that $m \geq m(\lambda)$, 
$\lambda \geq \sum_i \varepsilon_i a_{0}(h)/a_1(h)$.
Then by Lemma \ref{lem:lambda} and $(\flat 1,2)$, 
all generalized parabolic semi-stable sheaves 
$E_*$ are parametrized by $\Gamma^{ss}$,
$H^0(X,E(m)) \to H^0(X,E_i(m))$ is surjective and $H^1(X,E_i(m))=0$. 
For a sequence of quotients 
\begin{equation}
 V_m \otimes {\cal O}_X \to E(m) \to E_l(m) \to E_{l-1}(m) \to \dots
 \to E_1(m) \to 0 \in \Gamma^{ss},
\end{equation}
quotient vector spaces 
$\alpha:V_m \otimes W \to H^0(X,E(m+n))$ and 
$\alpha_i:V_m \to H^0(X,E_i(m))$ define a point of $G(n) \times \prod_i G_i$.
Thus we get a morphism $\Gamma^{ss} \to G(n) \times \prod_i G_i$.
As in \cite{Inaba:1}, we can show that this morphism is an immersion
\begin{equation}
 \Gamma^{ss} \hookrightarrow G(n) \times \prod_i G_i.
\end{equation}
$SL(V_m)$ acts on $G(n) \times \prod_i G_i$. 
Let ${\cal O}_{G(n)}(1)$ and ${\cal O}_{G_i}(1)$, $1 \leq i \leq l$ be
tautological line bundles on $G(n)$ and $G_i$ respectively.
These line bundles have $SL(V_m)$-linearizations.
We consider GIT semi-stability with respect to a ${\Bbb Q}$-line bundle
$L={\cal O}_{G(n)}(\beta_0) \otimes {\cal O}_{G_1}(\beta_1) \otimes
\dots \otimes {\cal O}_{G_l}(\beta_l)$.
\begin{prop}\cite[Prop. 3.2]{Inaba:1}
Let $\alpha:V_m \otimes W \to A$ and $\alpha_i:V_m \to A_i$
be quotients corresponding to a point of $G(n) \times \prod_i G_i$.
Then it is GIT semi-stable with respect to $L$
if and only if 
\begin{equation}
 \dim V_m (\beta_0 \dim \alpha(V' \otimes W)+\sum_i \beta_i 
 \dim \alpha_i(V'))-
 \dim V' (\beta_0 \dim \alpha(V_m \otimes W)+\sum_i \beta_i 
 \dim \alpha_i(V_m)) \geq 0
\end{equation}
for all non-zero subspaces $V'$ of $V_m$.
\end{prop}

We set $\beta_0:=(h(m)-\sum_{i=1}^l \varepsilon_i h_i(m))/a_1(h)n$,
$\beta_i:=\varepsilon_i$ for
$1 \leq i \leq l$.
We also set $V_i:=\ker(\alpha_{i|V'})$. Then 
\begin{equation}\label{eq:a1}
\begin{split}
&\dim V_m (\beta_0 \dim \alpha(V' \otimes W)+\sum_i \beta_i 
\dim \alpha_i(V'))-
\dim V' (\beta_0 \dim \alpha(V_m \otimes W)+\sum_i \beta_i 
\dim \alpha_i(V_m))\\
=&
h(m)\left(\frac{h(m)-\sum_i \varepsilon_i h_i(m)}{a_1(h)n}
\dim \alpha(V'\otimes W)+\sum_i\varepsilon_i
(\dim V'-\dim V_i)\right)\\
&-\dim V'\left(\frac{h(m)-\sum_i \varepsilon_i h_i(m)}{a_1(h)n}
h[m](n)+\sum_i \varepsilon_i h_i(m)\right)\\
=&
h(m)\left(\frac{h(m)-\sum_i \varepsilon_i h_i(m)}{a_1(h)n}
\dim \alpha(V'\otimes W)+\sum_i\varepsilon_i
(\dim V'-\dim V_i)\right.\\
&\left. -\dim V'-
\frac{(h(m)-\sum_i \varepsilon_i h_i(m))\dim V'}{a_1(h)n} \right)\\
= & h(m)\left(\frac{h(m)-\sum_i \varepsilon_i h_i(m)}{a_1(h)n}
\dim \alpha(V'\otimes W)-\sum_i\varepsilon_i \dim V_i-
\alpha_1 \dim V' \right.\\
& -\left.\frac{(h(m)-\sum_i \varepsilon_i h_i(m))\dim V'}{a_1(h)n} \right).
\end{split}
\end{equation}

Let $(G(n) \times \prod_i G_i)^{ss}$ be the open subscheme of
$G(n) \times \prod_i G_i$ consisting of GIT semi-stable points.
Then we get the following.
\begin{prop}
We set $\beta_0:=(h(m)-\sum_{i=1}^l \varepsilon_i h_i(m))/a_1(h)n$,
$\beta_i:=\varepsilon_i$ for
$1 \leq i \leq l$.
Then there is an integer $m_1$ such that for all $m \geq m_1$,
$\Gamma^{ss}$ is contained in $(G(n) \times \prod_i G_i)^{ss}$,
where $n \gg m$.
\end{prop}

\begin{proof}
We set 
\begin{equation}
{\cal F}:=\{E' \subset E(m)|\text{$E'$ is generated by $V' \subset V$ }\}.
\end{equation}
Since ${\cal F}$ is a bounded set, for a sufficiently large $n$ 
which depends on $m$,
\begin{enumerate}
\item
$\alpha(V' \otimes W)=H^0(X,E'(n))$,
$H^i(X,E'(n))=0$, $i>0$ and
\item
\begin{equation}
 \left|\frac{h(m)-\sum_i \varepsilon_i h_i(m)}{a_1(h)n}(\chi(E')-\dim V')
 \right|<\frac{1}{a_1(h)}.
\end{equation}
\end{enumerate}
Then $\dim \alpha(V' \otimes W)=\chi(E'(n))=a_1(E')n+\chi(E')$.
Since $E_*$ is semi-stable,
in the same way as in \cite[Prop. 2.5]{M-Y:1}, we see that
there is an integer $m_0$ such that for $m \geq m_0$ and
a generalized parabolic subsheaf $E_*'$ of
$E_*(m)$,
\begin{equation}\label{eq:pgs}
 \frac{\alpha_1 h^0(F_{l+1}(E'))+
 \sum_{i=1}^l \varepsilon_i h^0(F_i(E'))}{a_1(E')} \leq 
 \frac{\alpha_1 h^0(F_{l+1}(E(m)))+
 \sum_{i=1}^l \varepsilon_i h^0(F_i(E(m)))}{a_1(E)}
\end{equation}
and the equality holds, if and only if 
$\text{par-}\chi(E_*')/a_1(E')=\text{par-}\chi(E_*)/a_1(E)$.
Hence if the equality holds, then
$E'_*$ is semi-stable and we may assume that ($\flat1,2$) holds
for $E'_*(-m)$.
In particular, $\dim V'=\chi(E')$.  
We note that
\begin{equation}\label{eq:a6}
 \begin{split}
 &\frac{h(m)-\sum_i \varepsilon_i h_i(m)}{a_1(h)n}
 \dim \alpha(V'\otimes W)-\sum_i\varepsilon_i \dim V_i-
 \alpha_1 \dim V'-\frac{(h(m)-\sum_i \varepsilon_i h_i(m))\dim V'}{a_1(h)n}\\
 =& (h(m)-\sum_i \varepsilon_i h_i(m))\frac{a_1(E')}{a_1(E)}-
 (\sum_i\varepsilon_i \dim V_i+\alpha_1 \dim V')
 +\frac{h(m)-\sum_i \varepsilon_i h_i(m)}{a_1(h)n}(\chi(E')-\dim V').
 \end{split}
\end{equation}
By \eqref{eq:a1} and \eqref{eq:a6}, 
if the inequality in \eqref{eq:pgs} is strict, then
we get 
\begin{equation}\label{eq:a7}
\dim V_m (\beta_0 \dim \alpha(V' \otimes W)+\sum_i \beta_i 
\dim \alpha_i(V'))-
\dim V' (\beta_0 \dim \alpha(V_m \otimes W)+\sum_i \beta_i 
\dim \alpha_i(V_m)) > 0.
\end{equation}
If the equality holds in \eqref{eq:pgs}, then $\dim V'=\chi(E')$, and hence
\eqref{eq:a6} implies that L.H.S. of \eqref{eq:a7} is $0$.
Therefore our claim holds. 
\end{proof}

\begin{rem}
If $\beta_0=(h(m)-\sum_i \varepsilon_i h_i(m))/(a_1(h)n+t)$, 
$t \in {\Bbb Q}$, then we can show that
$\Gamma^{ss} \cap (G(n) \times \prod_i G_i)^{ss}$
parametrizes generalized semi-stable parabolic sheaves $E_*$ such that
\begin{equation}
 t\frac{\chi(E)}{a_1(E)} \leq t\frac{\chi(E')}{a_1(E')}
\end{equation}
for all generalized parabolic subsheaves $E'_*$ of 
$\text{par-}\chi(E'_*)/a_1(E')=\text{par-}\chi(E_*)/a_1(E)$.
If $t=h(m)$, then $L$ is nothing but the polarization in
\cite{Inaba:1}.
\end{rem}

\begin{prop}\label{prop:proper}
There is an integer $m_2$ such that for all $m \geq m_2$,
$\Gamma^{ss}$ is a closed subscheme of 
$(G(n) \times \prod_i G_i)^{ss}$,
where $n \gg m$.
\end{prop}

\begin{proof}
We choose an $m$ so that $h(m)/a_1(h)-\sum_i \varepsilon_i h_i(m)>0$.
We shall prove that 
$\Gamma^{ss} \to (G(n) \times \prod_i G_i)^{ss}$ is proper.
Let $(R,{\frak m})$ be a discrete valuation ring and $K$ the quotient field of $R$.
We set $T:=\Spec(R)$ and $U:=\Spec(K)$.
Let $U \to \Gamma^{ss}$ be a morphism such that $U \to \Gamma^{ss} \to 
(G(n) \times \prod_i G_i)^{ss}$ is extended to
a morphism $T \to (G(n) \times \prod_i G_i)^{ss}$.
Since $\Gamma$ is a closed subscheme of $\prod_i {\cal Q}_i$,
there is a morphism $T \to \Gamma$, {\it i.e},
there is a flat family of a sequence of quotients
\begin{equation}
 V \otimes {\cal O}_{T \times X} \to 
 {\cal E}(m) \to {\cal E}_l(m) \to {\cal E}_{l-1}(m) \to \cdots 
 \to {\cal E}_1(m) \to 0.
\end{equation}
Let $\alpha:V_m \otimes W \otimes R \to p_{T*}({\cal E}(m+n))$ and
$\alpha_i:V_m \otimes R \to A_i$ be the quotient bundles of
$V_m \otimes W \otimes R$ and $V_m \otimes R$ corresponding to
the morphism $T \to (G(n) \times \prod_i G_i)^{ss}$.
We set $E:={\cal E} \otimes R/{\frak m}$,
$E_i:={\cal E}_i \otimes R/{\frak m}$ and 
$F_{i+1}(E):=\ker(E \to E_i)$.
\begin{claim}
$V_m \to H^0(X,E(m))$ is injective.
\end{claim}
Indeed, we set $V':=\ker(V_m \to H^0(X,E(m)))$.
Then $\alpha(V' \otimes W)=0$.
By \eqref{eq:a1},

\begin{equation}
 \begin{split}
  0 \leq & \dim V_m (\beta_0 \dim \alpha(V' \otimes W)+\sum_i \beta_i 
  \dim \alpha_i(V'))-
  \dim V' (\beta_0 \dim \alpha(V_m \otimes W)+\sum_i \beta_i 
  \dim \alpha_i(V_m))\\
  =& h(m)\left(-\sum_i\varepsilon_i \dim V_i-
  \alpha_1 \dim V'-
  \frac{(h(m)-\sum_i \varepsilon_i h_i(m))\dim V'}{a_1(h)n}\right) \\
  \leq &-h(m)\frac{(h(m)-\sum_i \varepsilon_i h_i(m))\dim V'}{a_1(h)n}.
 \end{split}
\end{equation}
Therefore $V'=0$.

\begin{claim}\label{claim:2}
There is a rational number $\lambda$ which depends on $h_*$ and
$\alpha_*$ such that $E$ is of type $\lambda$.
\end{claim}
Proof of the claim:
Let $E \to E''$ be a quotient of $E$.
Let $E'$ be the kernel of $E \to E''$.
We note that $V_m \to H^0(X,E(m))$ is injective. 
We set $V':=V_m \cap H^0(X,E'(m))$.
Then $h^0(X,E''(m)) \geq \dim V_m -\dim V'$.
By \eqref{eq:a1}, 
\begin{equation}
 \begin{split}
  0 \leq & \dim V (\beta_0 \dim \alpha(V' \otimes W)+\sum_i \beta_i 
  \dim \alpha_i(V'))-
  \dim V' (\beta_0 \dim \alpha(V \otimes W)+\sum_i \beta_i 
  \dim \alpha_i(V))\\
  =& h(m)\left(\frac{h(m)-\sum_i \varepsilon_i h_i(m)}{a_1(h)n}
  \dim \alpha(V'\otimes W)+\sum_i\varepsilon_i
  (\dim V'-\dim V_i)-\dim V' \right.\\
  & -\left.\frac{(h(m)-\sum_i \varepsilon_i h_i(m))\dim V'}{a_1(h)n} \right).
\end{split}
\end{equation}

Let $F$ be a subsheaf of $E(m)$ generated by $V'$.
Then $F$ belongs to ${\cal F}$.  
Let $\varepsilon$ be a positive number such that
$h(m)/a_1(h)-\sum_i \varepsilon_i h_i(m)- \varepsilon>0$.
Since ${\cal F}$ is a bounded set, for a sufficiently large $n$ 
which depends on $m$ and $\varepsilon$,
we have $\alpha(V' \otimes W)=H^0(X,F(n))$,
$H^1(X,F(n))=0$ and 
\begin{equation}
 \left|\left(h(m)-\sum_i \varepsilon_i h_i(m) \right)
 \frac{\dim \alpha(V' \otimes W)}{a_1(h)n}
  -\left(h(m)-\sum_i \varepsilon_i h_i(m) \right)
 \frac{a_1(F)}{a_1(h)} \right|<
 \varepsilon.
\end{equation}
Therefore
\begin{equation}\label{eq:a2}
 \begin{split}
  0 \leq & \frac{h(m)-\sum_i \varepsilon_i h_i(m)}{a_1(h)n}
  \dim \alpha(V'\otimes W)+\sum_i\varepsilon_i
  (\dim V'-\dim V_i)-\dim V'\\
  < & (h(m)-\sum_i \varepsilon_i h_i(m))\frac{a_1(F)}{a_1(h)}+\varepsilon
  +\sum_i\varepsilon_i(\dim V'-\dim V_i)-\dim V',
 \end{split}
\end{equation}
where $V_i=\ker(\alpha_{i|V'})$.
Since $\dim(V'/V_i) \leq \dim (\im \alpha_i)=h_i(m)$,
we get 
\begin{equation}
 h(m)\frac{a_1(F)}{a_1(h)}-\dim V'>
 \sum_i \varepsilon_i h_i(m)\left(\frac{a_1(F)}{a_1(h)}-1 \right)-\varepsilon
 \geq -\sum_i \varepsilon_i h_i(m)-\varepsilon.
\end{equation}
Since $a_1(E') \geq a_1(F)$,
\begin{equation}
 \begin{split}
   \frac{h^0(X,E''(m))}{a_1(E'')} & \geq
   \frac{\dim V_m-\dim V'}{a_1(E'')}\\
  & >\left(\dim V_m-h(m)\frac{a_1(F)}{a_1(h)}-
  \sum_i \varepsilon_i h_i(m)-\varepsilon \right)
   \frac{1}{a_1(E'')}\\
  &=h(m)\frac{a_1(h)-a_1(F)}{a_1(h)}\frac{1}{a_1(E'')}-
   (\sum_i \varepsilon_i h_i(m)+\varepsilon)/a_1(E'')\\
  & \geq h(m)\frac{a_1(h)-a_1(E')}{a_1(h)}\frac{1}{a_1(E'')}-
   \sum_i \varepsilon_i h_i(m)-\varepsilon\\    
  &=\frac{h(m)}{a_1(h)}-\sum_i \varepsilon_i h_i(m)-\varepsilon>0.
 \end{split}
\end{equation}
There is a rational number $\lambda_1$ 
and an integer $m_0 \geq \lambda_1-a_{0}(h)/a_1(h)$ which depend on
$h(x)$, $\sum_i \varepsilon_i h_i(x)$ and $\varepsilon$ such that
\begin{equation}\label{eq:a5}
 \frac{h(m)}{a_1(h)}-\sum_i \varepsilon_i h_i(m) -\varepsilon \geq 
 \left(m+\frac{a_{0}(h)}{a_1(h)}-\lambda_1 \right)
\end{equation}
 for $m \geq m_0$.

By \cite[Lem. 1.17]{S:1}, there is a purely $1$-dimensional sheaf $G$ of
Hilbert polynomial $h(x)$
and a map $E \to G$ whose kernel is a coherent sheaf of dimension $0$.
Let $G \to G''$ be a quotient such that $G''$ is semi-stable.
We set $E':=\ker(E \to G'')$ and $E'':=\im(E \to G'')$.
Since $h^0(X,E''(m)) \leq h^0(X,G''(m))$ and $a_1(E'')=a_1(G'')$,
\begin{equation}
 \begin{split}
  \frac{h^0(X,G''(m))}{a_1(G'')} & \geq \frac{h^0(X,E''(m))}{a_1(E'')} \\    
  & \geq \frac{h(m)}{a_1(h)}-\sum_i \varepsilon_i h_i(m)-\varepsilon>0.
 \end{split}
\end{equation}
Since $G''$ is semi-stable, 
\cite[Cor. 1.7]{S:1} implies that 
\begin{equation}
 \frac{h^0(X,G''(m))}{a_1(G'')} \leq 
  \begin{cases}
   0, & c+m+\frac{\chi(G'')}{a_1(G'')}<0\\
   c+m+\frac{\chi(G'')}{a_1(G'')},
   & c+m+\frac{\chi(G'')}{a_1(G'')} \geq 0,
  \end{cases}
\end{equation}
where $c$ is a constant which only depends on $a_1(h)$.
Since $h^0(X,G''(m))>0$, we get $c+m+\chi(G'')/a_1(G'') \geq 0$ and
$\chi(G'')/a_1(G'') \geq a_0(h)/a_1(h)-\lambda_1-c$,
which means that $G$ is of type $\lambda:=a_0(h)/a_1(h)-\lambda_1-c$.
Replacing $m$, we may assume that for all type $\lambda$ sheaves $I$
of Hilbert polynomial $h(x)$, 
$I(m)$ is generated by global sections and $H^i(X,I(m))=0$, $i>0$.
In particular $h^0(X,G(m))=h(m)=\dim V_m$.
Assume that $H^0(X,E(m)) \to H^0(X,G(m))$ is not injective and let
$V'$ be the kernel.
Then we get a contradiction from the inequality \eqref{eq:a2} 
and $\alpha_1>0$.
Thus $H^0(X,E(m)) \to H^0(X,G(m))$ is injective, and hence it is
isomorphic.
Since $G(m)$ is generated by global sections,
$E \to G$ must be surjective, which implies that it is isomorphic.
Therefore $E$ is of pure dimension $d$, of type $\lambda$ and 
$V_m \to H^0(X,E(m))$ is an isomorphism.
Thus we complete the proof of Claim \ref{claim:2}. 

We assume that $m \geq m(\lambda)$.
Then $H^0(X,E(m)) \to H^0(X,E_i(m))$ is surjective.
Thus $\beta_i:V \otimes R \to p_{T*}({\cal E}_i(m))$ is surjective
and define a morphism $T \to G_i$.
Since $\beta_{i|U}=\alpha_{i|U}$ as elements of $G_i$, 
we get $\beta_i=\alpha_i$.
Assume that there is a generalized parabolic quotient $E_* \to E''_*$
which destabilizes semi-stability.
Since $E_*':=\ker (E_* \to E_*'')$
satisfies $(\flat 1, 2)$, we get that
$V'=H^0(X,E'(m))$, $\alpha_i(V')=H^0(X,E_i'(m))$ and  
\begin{equation}
 \frac{\text{par-}\chi(E_*''(m))}{a_1(E'')} \geq 
 \frac{\text{par-}\chi(E_*(m))}{a_1(h)}-\varepsilon.
\end{equation}
Since $\varepsilon_i$ are rational numbers,
for a sufficiently small $\varepsilon$,
we get  
\begin{equation}\label{eq:ss}
\frac{\text{par-}\chi(E_*''(m))}{a_1(E'')} \geq 
\frac{\text{par-}\chi(E_*(m))}{a_1(h)}, 
\end{equation}
which is a contradiction.
Therefore $E$ is generalized parabolic semi-stable.
Thus we get a lifting of $T \to \Gamma^{ss}$ and 
conclude that $\Gamma^{ss} \to (G(n) \times \prod_i G_i)^{ss}$ is proper.
\end{proof}
By standard arguments, we see that 
$SL(V_m)s, s \in \Gamma^{ss}$ is a closed orbit
if and only if 
the corresponding generalized parabolic semi-stable sheaf $(E,F(E),\alpha_*)$
is isomorphic to
$\oplus_i(E_i,F(E_i), \alpha_*)$,
where $(E_i,F(E_i),\alpha)$ are generalized parabolic stable sheaves.

\begin{thm}\label{thm:parabolic}
There is a moduli scheme $\overline{M}_{X/{\Bbb C}}^{h_*,\alpha_*}$
parametrizing $S$-equivalence classes of
generalized parabolic semi-stable sheaves $(E,F(E),\alpha_*)$ of
Hilbert polynomials $h_*$.
\end{thm} 

Let $D$ be a Cartier divisor on $X$.
Then 
\begin{equation}
\{(E,F(E),\alpha_*)| \text {$E(-D) \to E \to E_l$ is
a 0 map} \}
\end{equation}
has a natural closed subscheme structure of 
$\overline{M}_{X/{\Bbb C}}^{h_*,\alpha_*}$.
Assume that 
$h_l(m)=\chi(E(m))-\chi(E(-D)(m))$.
Then if $\dim(D \cap \Supp(E))=0$, then $E(-D) \to E$ is injective
and the image is $F_{l+1}(E)$.
Thus $(E,F(E),\alpha_*)$ becomes a parabolic semi-stable sheaf.  
Since this condition is an open condition, 
we get $\overline{M}_{D/X/{\Bbb C}}^{h_*,\alpha_*}$
as a locally closed subscheme of $\overline{M}_{X/{\Bbb C}}^{h_*,\alpha_*}$.
\begin{thm}
Let $D$ be a Cartier divisor on $X$ and assume that 
$h_l(m)=\chi(E(m))-\chi(E(-D)(m))$.
Then there is a moduli scheme $\overline{M}_{D/X/{\Bbb C}}^{h_*,\alpha_*}$
parametrizing $S$-equicalence classes of
parabolic semi-stable sheaves $(E,F(E),\alpha_*)$ of
Hilbert polynomials $h_*$.
\end{thm}

Assume that $X$ is a surface.
Let $\gamma$ be an element of $K(X)$ such that $\rk \gamma=0$
and $c_1(\gamma)$ is effective.
We set $h(x)=(c_1(\gamma),{\cal O}_X(1))x+\chi(\gamma)$.
Let $D$ is a irreducible and reduced 
divisor such that $(D,{\cal O}_X(1))>(c_1(\gamma),{\cal O}_X(1))$.
Then $\overline{M}_{D/X/{\Bbb C}}^{h_*,\alpha_*}$ becomes compact.
Therefore we get Theorem \ref{thm:alpha-ss}.

\begin{rem}
Although we assume that $X$ is defined over a field $k$,
we can easily generalize our construction to relative setting
$X \to S$, where $S$ is of finite type over $k$ and 
$X \to S$ is projective.
\end{rem}

\vspace{1pc}

{\it Acknowledgement.}
I would like to thank M. Inaba for explaining the difficulty
of constructing the moduli space of parabolic semi-stable sheaves.

\end{document}